\input ./style/arxiv-general.cfg
\documentclass[MSNbibl,number,citesort,seceqn,dvips]{arxbj}
\makeatletter
   \@ifpackageloaded{graphicx}{}{\usepackage{graphicx}}
\makeatother

\usepackage[dvips,all]{xy}
\usepackage{xypic}
\usepackage{graphicx}

\usepackage{mathrsfs}

\aid{0}
\volume{21}
\issue{3}
\pubyear{2015}
\firstpage{1386}
\lastpage{1411}
\doi{10.3150/14-BEJ607} % kopijuoti is 'New paper accepted'
\makeatletter

\newcommand{\eqref}[1]{(\ref{#1})}
\renewcommand{\emptyset}{\varnothing}

\newremark{rmk}{Remark}[section]
\newremark{example}{Example}[section]
\newtheorem{lemma}{Lemma}[section]
\newtheorem{prop}{Proposition}[section]
\newtheorem{cor}{Corollary}[section]
\newtheorem{them}{Theorem}[section]
\newproclaim{defn}{Definition}[section]

\newcommand{\idN}{\mathbf{Id}}
\newcommand{\Dmn}{{\mathbf{D}_{m,n}}}
\newcommand{\Lip}{\operatorname{Lip}}
\newcommand{\bM}{\mathbf{M}}
\newcommand{\Frag}{\operatorname{Frag}}
\newcommand{\Coag}{\operatorname{Coag}}
\newcommand{\partitionsNk}{{\mathcal{P}_{\mathbb{N}:k}}}
\newcommand{\partitionsN}{{\mathcal{P}_{\mathbb{N}}}}
\newcommand{\matrixNk}{{\mathcal{M}_{\mathbb{N}:k}}}
\newcommand{\mathL}{{\mathcal{L}_{\mathbb{N}:k}}}
\newcommand{\diag}{\operatorname{diag}}
\newcommand{\symmetricn}{{\mathscr{S}_n}}
\newcommand{\symmetrick}{{\mathscr{S}_k}}
\newcommand{\stochk}{{\mathcal{S}_k}}
\newcommand{\simplexk}{{\Delta_k}}
\newcommand{\equalinlaw}{=_{\mathcal{L}}}
\makeatother

\begin{document}
\begin{frontmatter}

\title{Lipschitz partition processes}
\runtitle{Lipschitz partition processes}

\begin{aug}
%%%% inicialai - be tarpu
\author{\inits{H.}\fnms{Harry}~\snm{Crane}\corref{}\ead[label=e1]{hcrane@stat.rutgers.edu}}% \and
%\author{\inits{}\fnms{}~\snm{}\thanksref{}\ead[label=e2]{}}
%\author{\inits{}\fnms{}~\snm{}}
%%\runauthor{} %% auto
%\dedicated{}
\address{Department of Statistics \&
Biostatistics, Rutgers University, 110 Frelinghuysen Road, Hill Center, Room 461,
Piscataway, NJ 08854, USA. \printead{e1}}
%\address[]{}
\end{aug}

% HISTORY:
\received{\smonth{9} \syear{2012}}
\revised{\smonth{9} \syear{2013}}

% ABSTRACT
%
\begin{abstract}
We introduce a family of Markov processes on set partitions with a
bounded number of blocks, called \emph{Lipschitz partition processes}.
We construct these processes explicitly by a Poisson point process on
the space of Lipschitz continuous maps on partitions. By this
construction, the \emph{Markovian consistency property} is readily
satisfied; that is, the finite restrictions of any Lipschitz partition
process comprise a compatible collection of finite state space Markov
chains. We further characterize the class of exchangeable Lipschitz
partition processes by a novel set-valued matrix operation.
\end{abstract}

% KEYWORDS
% visi is mazosios raides ir pagal abecele
%

\end{frontmatter}

%s1 #&#
\section{Introduction}\label{section:introduction}

Partition-valued Markov processes, particularly coalescent and
fragmentation processes, arise as mathematical models in population
genetics and mathematical biology. Initially, Ewens \cite{Ewens1972}
derived his celebrated sampling formula while studying neutral allele
sampling in population genetics.
Extending Ewens's work, Kingman characterized exchangeable partitions
of the natural numbers \cite{Kingman1978a,Kingman1978b}, which play a
larger role in the mathematical study of genetic diversity \cite
{Kingman1980}. Related applications in phylogenetics and the study of
ancestral lineages prompted Kingman's coalescent process \cite
{Kingman1982}, which arises as the scaling limit of both Wright--Fisher
and Moran models under different regimes \cite{TavareStFlour}.
Exchangeable coalescent and fragmentation processes have also taken
hold in the probability literature because of some beautiful
relationships to classical stochastic process theory, for example,
Brownian motion and L\'evy processes. For specific content in the
literature, see
\cite
{AldousPitmanAdditive,Berestycki2004,Bertoin2001a,Bertoin2002,Pitman1995};
for recent
overviews of this theory, see \cite{Bertoin2006,Pitman2005}.

In this paper, we study a family of Markov processes on labeled
partitions with a finite number $k\geq1$ of classes. By a simple
projection, we describe a broad class of processes on the space of
partitions with at most $k$ blocks.
Processes on this space are cursorily related to composition structures
for ordered partitions, for example, \cite
{DonnellyJoyce1991,Gnedin1997}, but our approach more closely follows
previous work \cite{Crane2011a}, which is motivated by DNA sequencing
applications.
In addition to genetics applications, processes on this subspace relate
to problems of cluster detection and classification in which the total
number of classes is finite, for example,
\cite{BoothCasellaHobert2008,Crane2014,McCullaghYang2008}.

Prior to \cite{Crane2011a}, coagulation--fragmentation processes
dominated the literature.
The processes in \cite{Crane2011a} do not evolve by fragmentation or
coagulation; their jumps involve simultaneous fragmentation and
coagulation of all blocks.
To describe a broader class of processes, we incorporate ideas from the
coagulation--fragmentation literature as well as our previous work.
We call these \emph{Lipschitz partition processes}.

Our main theorems are \emph{not} corollaries of the many results for
fragmentation and coalescent processes.
Instead, our approach extracts fundamental properties of these
processes, specifically their construction from the $\Coag$ and $\Frag
$ operators; see, for example, Bertoin \cite{Bertoin2006}, Chapters~3--4.
Importantly, these operators are Lipschitz continuous and
associative. From these observations, we construct a family of
processes by repeated application of random Lipschitz continuous maps
that act on the space of partitions.

In the exchangeable case, the random maps are confined to the subspace
of \emph{strongly Lipschitz continuous} functions, which we
characterize in full by a class of specially structured set-valued
matrices (Section~\ref{section:matrix multiplication}).
These set-valued matrices act on labeled partitions similarly to the
action of a matrix on a real-valued vector (with obvious modifications
to the operations addition and multiplication).
They also establish an intimate connection between exchangeable
Lipschitz partition processes and random stochastic matrices
(Section~\ref{section:measure-valued}).

%s1.1 #&#
\subsection{General construction: Overview}

For now, we regard a labeled partition as a finite collection of
non-overlapping, labeled subsets.

Consider the following construction of a discrete-time Markov chain.
Let $\Lambda_0$ be an initial state and let $F_1,F_2,\ldots$ be
independent and identically distributed (i.i.d.) random maps on the
space of labeled partitions. Then, for each $t\geq1$, we define
%
%e1.1 #&#
%
\begin{equation}
\label{eq:general-discrete}\Lambda_t:=F_t(\Lambda _{t-1})=(F_t
\circ F_{t-1}\circ\cdots\circ F_1) (\Lambda_0).
\end{equation}
The collection $\boldsymbol{\Lambda}:=(\Lambda_t, t\geq0)$ is a
discrete-time Markov chain.

We study an analogous construction for continuous-time processes.
Instead of an i.i.d. sequence of random maps, we construct $\boldsymbol
{\Lambda}$ from a Poisson point process on the space of maps.
Informally, if $\mathbf{F}:=\{(t,F_t)\}$ is a realization of such a
Poisson point process (where each $F_t$ is a map), we construct
$\boldsymbol{\Lambda}$ by putting
%
%e1.2 #&#
%
\begin{equation}
\label{eq:general-PPP} \Lambda_t:=\cases{ %
%\begin{array}{cc}
F_t(\Lambda_{t-}),&\quad$t\mbox{ is an atom time of }
\mathbf{F}$,
\cr
\Lambda_{t-},& \quad$\mbox{otherwise}$, %\end{array}
%\right.
%,
}\qquad\mbox{for every }t>0.
\end{equation}

We are interested in processes $\boldsymbol{\Lambda}$ that exhibit
\begin{itemize}
\item\emph{Markovian consistency}, that is, for each $n\in\mathbb
{N}$, the restriction of $\boldsymbol{\Lambda}$ to labeled partitions of
$[n]:=\{1,\ldots,n\}$ is a Markov chain, and
\item\emph{exchangeability}, that is, the law of $\boldsymbol
{\Lambda}$
is invariant under relabeling of elements of $\mathbb{N}$.
\end{itemize}
Markovian consistency might also be called the \emph{projective Markov
property}, meaning the projection of $\boldsymbol{\Lambda}$ to spaces of
finite labeled partitions is also Markov. Throughout the paper, we use
the term \emph{consistency} in place of \emph{Markovian consistency}.
Consistency plays a central role not only in this paper but also more
widely in the study of partition-valued Markov processes. In general, a
function of a Markov process need not be Markov, and so consistency is
not trivially satisfied; see Example~\ref{ex:counterexample}.

We pay special attention to the exchangeable case, for which we can
make some precise statements. In this case, we show that the Poisson
point process $\mathbf{F}$ is supported on the space of maps having
the \emph{strong Lipschitz property} (Section~\ref{section:strong lipschitz}).

The general approach outlined in \eqref{eq:general-discrete} and
\eqref{eq:general-PPP} can be applied to construct processes on the
unrestricted space of set partitions, or even ordered set partitions,
but we do not treat these cases.
In our main theorems, we show a correspondence between strongly
Lipschitz maps on labeled partitions and $k\times k$ set-valued matrices.
Without bounding the number of classes, we cannot obtain such a precise
statement.

%s1.2 #&#
\subsection{Organization of the paper}

We organize the paper as follows. In Section~\ref{section:preliminaries},
we give some preliminaries for partitions and
labeled $k$-partitions. In Section~\ref{section:ct Markov}, we
introduce the general class of Lipschitz partition processes; and in
Section~\ref{section:exchangeable processes}, we specialize to
exchangeable Lipschitz partition processes. In
Section~\ref{section:discrete-time processes}, we discuss
discrete-time Markov
chains. In Section~\ref{section:concluding remarks}, we make some
concluding remarks about projections to unlabeled set partitions and
more general issues concerning partition-valued Markov processes.

%s2 #&#
\section{Preliminaries}\label{section:preliminaries}

%s2.1 #&#
\subsection{Partitions}\label{section:partitions}

For $n\in\mathbb{N}=\{1,2,\ldots\}$, a \emph{partition} $\pi$ of
$[n]:=\{1,\ldots,n\}$ is a collection $\{b_1,\ldots,b_r\}$ of
non-empty, disjoint subsets (blocks) satisfying $\bigcup_{i=1}^r
b_i=[n]$. Alternatively, $\pi$ can be regarded as an equivalence
relation $\sim_\pi$, where
%
%e2.1 #&#
%
\begin{equation}
\label{eq:equiv relat}i\sim_\pi j \quad\Longleftrightarrow\quad i\mbox{ and }j
\mbox{ are in the same block of }\pi.
\end{equation}
We write $\#\pi$ to denote the number of blocks of $\pi$. Unless
otherwise stated, we assume that the blocks of $\pi$ are listed in
increasing order of their least element. We write ${\mathcal
{P}_{[n]}} $ to denote the space of partitions of $[n]$.

Writing $\symmetricn$ to denote the symmetric group acting on $[n]$,
we define the \emph{relabeling} $\pi\in{\mathcal
{P}_{[n]}} $ by $\sigma\in\symmetricn$, $\pi\mapsto\pi
^{\sigma}$, where
\[
i\sim_{\pi^{\sigma}}j \quad\Longleftrightarrow\quad\sigma(i)\sim
_\pi\sigma(j).
\]
Furthermore, for $m\leq n$, we define the \emph{restriction} of $\pi
\in{\mathcal{P}_{[n]}} $ to $\mathcal{P}_{[m]}$ by
\[
\pi_{|[m]}=\Dmn\pi:=\bigl\{b\cap[m]\dvtx b\in\pi\bigr\}\setminus\{
\emptyset\},
\]
the restriction of each block of $\pi$ to $[m]$ after removal of any
empty sets. In general, to any injective map $\psi\dvtx[m]\rightarrow
[n]$, we associate a projection $\psi'\dvtx{\mathcal
{P}_{[n]}} \rightarrow\mathcal{P}_{[m]}$, where
\[
i\sim_{\psi'(\pi)} j \quad\Longleftrightarrow\quad\psi(i)\sim _\pi
\psi(j).
\]

We write $\partitionsN$ to denote the space of partitions of $\mathbb
{N}$, which are defined as compatible sequences $(\pi_n, n\in\mathbb
{N})$ of finite set partitions.
For $m\leq n$, we say $\pi\in{\mathcal{P}_{[n]}} $
and $\pi'\in\mathcal{P}_{[m]}$ are \emph{compatible} if $\pi
_{|[m]}=\pi'$; and we call $(\pi_n, n\in\mathbb{N})$ a \emph
{compatible sequence} if $\pi_n\in{\mathcal
{P}_{[n]}} $ and $\pi_m=\Dmn\pi_n$, for all $m\leq n$, for
every $n\in\mathbb{N}$.

Writing $\mathbf{n}(\pi,\pi'):=\max\{n\in\mathbb{N}\dvtx \pi
_{|[n]}=\pi_{|[n]}'\}$, we equip $\partitionsN$ with the ultrametric
%
%e2.2 #&#
%
\begin{equation}
\label{eq:ultrametric partition}{d}_{\partitionsN
}\bigl(\pi,\pi'\bigr):=2^{-\mathbf{n}(\pi,\pi')},
\qquad\pi,\pi'\in \partitionsN,
\end{equation}
under which $(\partitionsN,d_{\partitionsN})$ is complete, separable,
and naturally endowed with the discrete $\sigma$-field $\sigma
\langle\bigcup_{n\in\mathbb{N}}{\mathcal
{P}_{[n]}}  \rangle$.

%s2.2 #&#
\subsection{Random partitions}\label{section:random partitions}

A sequence $(\mu_n, n\in\mathbb{N})$ of measures on the system
$({\mathcal{P}_{[n]}}, n\in\mathbb{N})$, where
$\mu_n$ is a measure on ${\mathcal{P}_{[n]}} $ for
each $n\in\mathbb{N}$, is \emph{consistent} if
%
%e2.3 #&#
%
\begin{equation}
\label{eq:consistent measure-partition} \mu_m=\mu_n\mathbf{D}_{m,n}^{-1}
\qquad\mbox{for every }m\leq n;
\end{equation}
that is, $\mu_m$ coincides with the law $\mu_n\mathbf{D}_{m,n}^{-1}$
induced by the restriction map. By Kolmogorov's extension theorem, any
consistent collection of measures determines a unique measure $\mu$ on~$\partitionsN$. This circle of ideas is central to the theory of
random partitions of $\mathbb{N}$ as it permits the explicit
construction of a random partition $\Pi$ through its compatible
sequence $(\Pi_n, n\in\mathbb{N})$ of finite random partitions.

A random partition $\Pi$ of $\mathbb{N}$ is called \emph
{exchangeable} if $\Pi^{\sigma}\equalinlaw\Pi$ for all permutations
$\sigma\dvtx\mathbb{N}\rightarrow\mathbb{N}$ that fix all but finitely
many elements of $\mathbb{N}$, where $\equalinlaw$ denotes \emph
{equality in law}. Kingman \cite{Kingman1978b} gives a de Finetti-type
characterization of exchangeable random partitions of $\mathbb{N}$
through the paintbox process. Let
\[
{\Delta^{\downarrow}}:= \Biggl\{(s_1,s_2,\ldots)\dvtx
s_1\geq s_2\geq\cdots\geq 0, \sum
_{i=1}^{\infty}s_i\leq1 \Biggr\}
\]
denote the space of \emph{ranked mass partitions}. Given $s\in
{\Delta^{\downarrow}}$, we write $s_0:=1-\sum_{i\geq1}s_i$ and construct
$\Pi$ as follows. Let $X_1,X_2,\ldots$ be a sequence of independent
(but not necessarily identically distributed) random variables with law
\[
P_s\{X_i=j\}:=\cases{ %
%\begin{array}{cc}
s_j,&\quad$j\geq1$,
\cr
s_0,&\quad$j=-i$,
\cr
0,&\quad$
\mbox{otherwise}$. %\end{array}
%%
%\right.
}
\]
Given $X:=(X_1,X_2,\ldots)$, we define $\Pi:=\Pi(X)$ by the relation
\[
i\sim_{\Pi} j \quad\Longleftrightarrow\quad X_i=X_j.
\]
We write $\varrho_{s}$ to denote the law of $\Pi$, called a \emph
{paintbox process directed by} $s$. For $n\in\mathbb{N}$, we write
$\varrho_s^{(n)}$ to denote the restriction of $\varrho_s$ to a
probability measure on ${\mathcal{P}_{[n]}} $. In
this way, $(\varrho_s^{(n)}, n\in\mathbb{N})$ is a consistent
collection of finite-dimensional measures determining $\varrho_s$.
More generally, given a measure $\nu$ on ${\Delta^{\downarrow}}$,
the $\nu
$-mixture of paintbox processes is defined by
\[
\varrho_{\nu}(\cdot):=\int_{{\Delta^{\downarrow}}}\varrho
_{s}(\cdot)\nu(\mathrm{d}s).
\]
Kingman's correspondence associates every exchangeable random partition
of $\mathbb{N}$ with a unique probability measure on ${\Delta
^{\downarrow}}$.

A widely circulated example of a sequential construction is the Chinese
restaurant process.
Overall, the Chinese restaurant process produces a compatible
collection $(\Pi_n, n\in\mathbb{N})$ of finite partitions for which
each $\Pi_n$ obeys the Ewens distribution on ${\mathcal
{P}_{[n]}} $.
The random partition $\Pi$ determined by $(\Pi_n, n\in\mathbb{N})$
obeys the \emph{Ewens process}, whose directing measure is the
two-parameter Poisson--Dirichlet distribution with parameter $(0,\theta
)$; see \cite{Pitman2005} for more information on the distinguishing
properties of the Ewens distribution.

%s2.3 #&#
\subsection{Partition-valued Markov processes}\label
{section:partition-valued process}

In this paper, we study Markov processes $\boldsymbol{\Pi}:=(\Pi_t,
t\geq0)$ on $\partitionsN$ that are
\begin{itemize}
\item\emph{consistent}: for each $n\in\mathbb{N}$, $\boldsymbol
{\Pi
}_{|[n]}:=(\Pi_{t|[n]}, t\geq0)$ is a Markov chain on $
{\mathcal{P}_{[n]}} $; and
\item\emph{exchangeable}: $\boldsymbol{\Pi}^{\sigma}:=(\Pi
^{\sigma
}_t, t\geq0)\equalinlaw\Pi$ for all permutations $\sigma\dvtx
\mathbb
{N}\rightarrow\mathbb{N}$ that fix all but finitely many $n\in
\mathbb{N}$.
\end{itemize}
In this case, exchangeability refers to joint exchangeability in the
sense that elements are relabeled according to the same partition at
all time points. Consistency refers to a preservation of the Markov property.

A consistent Markov process $\boldsymbol{\Pi}$ on $\partitionsN$ can be
constructed sequentially through its finite restrictions $(\boldsymbol
{\Pi
}_{|[n]}, n\in\mathbb{N})$, but care must be taken to ensure that
each of the restrictions $\boldsymbol{\Pi}_{|[n]}$ has c\`adl\`ag sample
paths. Perhaps the most well-known example of an exchangeable and
consistent Markov process on $\partitionsN$ is the exchangeable
coalescent process.

%s2.3.1 #&#
\subsubsection{Exchangeable coalescent process}\label{section:coalescent}

The construction of the coalescent process from the $\Coag$-operator
telegraphs our general approach.
Let $\pi:=\{b_1,b_2,\ldots\}$ be any partition of a finite or
infinite set with $\#\pi=k\in\mathbb{N}\cup\{\infty\}$, and let
$b':=\{b_1',b_2',\ldots\}$ be a partition of $[k']$, for any $k'\geq
k$. We call $\pi'':=\Coag(\pi,\pi'):=\{b_1'',b_2'',\ldots\}$ the
\emph{coagulation of $\pi$ by $\pi'$}, where
%
%e2.4 #&#
%
\begin{equation}
\label{eq:coag} b_i'':=\bigcup
_{j\in b_i'}b_j,\qquad i\geq1.
\end{equation}
(To maintain the definition of $\pi''$ as a partition, we remove any
empty sets that result from this operation.) Essential to definition
\eqref{eq:coag} is that the blocks of $\pi$ are ordered in ascending
order of their least element. For example, let $\pi:=1356/2/47/8$ and
$\pi':=135/24$, then
\[
\Coag\bigl(\pi,\pi'\bigr)=\Coag(1356/2/47/8,135/24)=134567/28.
\]
In words: block $\{1,3,5\}$ of $\pi'$ indicates that we merge the
\emph{first}, \emph{third}, and \emph{fifth} blocks of $\pi$, while
block $\{2,4\}$ indicates that we merge the \emph{second} and \emph
{fourth} blocks of $\pi$. (We ignore any elements of $\pi'$ larger
than $\#\pi$; for example, there is no fifth block of $\pi$ and so
the position of $5$ in $\pi'$ does not affect $\Coag(\pi,\pi')$.)

The $\Coag$ operator has been used extensively in the study of
coalescent processes; see Chapter~4 of Bertoin \cite{Bertoin2006}.
Let $\mu$ be a measure on $\partitionsN$ such that
%
%e2.5 #&#
%
\begin{equation}
\label{eq:regularity coalescent} \mu\bigl(\{\mathbf{0}_{\mathbb{N}}\}\bigr)=0 \quad\mbox{and}
\quad\mu\bigl(\{\pi \in \partitionsN\dvtx\pi_{|[n]}\neq
\mathbf{0}_{[n]}\}\bigr)<\infty\qquad \mbox{for every }n\in\mathbb{N},
\end{equation}
where $\mathbf{0}_{A}$ denotes the partition of $A\subseteq\mathbb
{N}$ into singletons.
Also, let $\mathbf{B}:=\{(t,B_t)\}\subset[0,\infty)\times
\partitionsN$ be a Poisson point process with intensity $\mathrm{d}t\otimes\mu
$ (where $\mathrm{d}t$ denotes Lebesgue measure on $[0,\infty)$). Given
$\mathbf{B}$, we construct a coalescent process $\boldsymbol{\Pi
}:=(\Pi
_t, t\geq0)$ on $\partitionsN$ as follows. For each $n\in\mathbb
{N}$, we specify $\boldsymbol{\Pi}^{[n]}:=(\Pi^{[n]}_t, t\geq0)$ on
${\mathcal{P}_{[n]}} $ by $\Pi^{[n]}_0=\mathbf
{0}_{[n]}$ and, for every $t>0$,
\begin{itemize}
\item if $t>0$ is an atom time of $\mathbf{B}$ such that
$B_{t|[n]}\neq\mathbf{0}_{[n]}$, then we put $\Pi^{[n]}_{t}=\Coag
(\Pi^{[n]}_{t-},B_t)$;
\item otherwise, we put $\Pi^{[n]}_{t}=\Pi^{[n]}_{t-}$.
\end{itemize}
Note that, by the definition of $\Coag$, $\Pi^{[m]}_t=\Dmn\Pi
^{[n]}_t$ for all $t\geq0$, for all $m\leq n$. Hence, $(\boldsymbol
{\Pi
}^{[n]}, n\in\mathbb{N})$ is a compatible collection of processes.
Furthermore, by \eqref{eq:regularity coalescent}, each $\boldsymbol
{\Pi
}^{[n]}$ is a Markov chain on ${\mathcal{P}_{[n]}} $
with c\`adl\`ag sample paths. Hence, $(\boldsymbol{\Pi}^{[n]}, n\in
\mathbb{N})$ determines a consistent Markov process $\boldsymbol{\Pi}$
on $\partitionsN$. If, in addition, $\mu$ is exchangeable, then
$\boldsymbol{\Pi}$ is exchangeable. The process constructed in this way
is called a \emph{coalescent process}.

%re2.1 #&#
%
\begin{rmk}
The construction of $\boldsymbol{\Pi}$ from the collection
$(\boldsymbol{\Pi
}^{[n]}, n\in\mathbb{N})$ of finite state space processes, rather
than directly from the entire process $\mathbf{B}$, is necessary. In
general, \eqref{eq:regularity coalescent} permits $\mathbf{B}$ to
have infinitely many atoms in arbitrarily small intervals of $[0,\infty
)$; but, by the second half of \eqref{eq:regularity coalescent}, there
can be only finitely many atom times $t>0$ for which $B_{t|[n]}\neq
\mathbf{0}_{[n]}$, for each $n\in\mathbb{N}$. Therefore, while the
Poisson point process construction cannot be applied directly to
construct $\boldsymbol{\Pi}$ (because the atom times might be dense in
$[0,\infty)$), we can construct $\boldsymbol{\Pi}$ sequentially by
building a compatible collection of processes that are consistent in
distribution.
\end{rmk}

An important property of the $\Coag$ operator is Lipschitz continuity
with respect to \eqref{eq:ultrametric partition}, that is, for every
$\pi
\in\partitionsN$,
\[
d_{\partitionsN}\bigl(\Coag\bigl(\pi',\pi\bigr),\Coag\bigl(
\pi'',\pi\bigr)\bigr)\leq d_{\partitionsN}\bigl(
\pi',\pi''\bigr) \qquad\mbox{for all }
\pi',\pi''\in \partitionsN.
\]
Furthermore, $\Coag\dvtx\partitionsN\times\partitionsN\rightarrow
\partitionsN$ is associative in the sense that
\[
\Coag\bigl(\pi,\Coag\bigl(\pi',\pi''
\bigr)\bigr)=\Coag\bigl(\Coag\bigl(\pi,\pi'\bigr),
\pi''\bigr) \qquad\mbox{for all }\pi,
\pi',\pi''\in\partitionsN.
\]
Lipschitz continuity is important for the consistency property because
it implies that the coagulation of $\pi_{|[n]}$ by $\pi'$ depends
only on $\pi_{|[n]}'$, for every $n\in\mathbb{N}$. Associativity
ensures the construction of $\boldsymbol{\Pi}$ is well-defined.

The $\Frag$ operator acts as the dual to $\Coag$ in the related study
of fragmentation processes. Analogously to the above construction, the
$\Frag$ operator can be used to construct fragmentation processes on
$\partitionsN$, but we do not discuss those details. We only
acknowledge that the $\Frag$ operator is also Lipschitz continuous
with respect to \eqref{eq:ultrametric partition}. These operators are
important because they characterize the semigroup of coagulation and
fragmentation processes. By Lipschitz continuity, the semigroups of
these processes are easily shown to fulfill the Feller property (under
the additional regularity condition \eqref{eq:regularity coalescent},
or its analog for fragmentation processes).

The coalescent process above need not be exchangeable. Bertoin \cite
{Bertoin2006} only considers the exchangeable case and so specializes
to the case in which $\mu$ in \eqref{eq:regularity coalescent} is the
directing measure of a paintbox process.

%s2.3.2 #&#
\subsubsection{\texorpdfstring{Processes on partitions with a bounded number of blocks
(Crane \cite{Crane2011a}, Section~4.1)}{Processes on partitions with a bounded number of blocks
(Crane [8], Section~4.1)}}\label{section:self-similar}

For $k\in\mathbb{N}$, let $\partitionsNk:=\{\pi\in\partitionsN
\dvtx\#
\pi\leq k\}$ be the subcollection of partitions of $\mathbb{N}$ with
$k$ or fewer blocks, and let ${\Delta_k^{\downarrow}}:= \{
(s_1,\ldots
,s_k)\dvtx s_1\geq\cdots\geq s_k\geq0, \sum_{i=1}^k s_i=1 \}$
denote the \emph{ranked $k$-simplex}. For any probability measure $\nu
$ on ${\Delta_k^{\downarrow}}$, the paintbox measure $\varrho_{\nu
}$ is
supported on $\partitionsNk$.

In \cite{Crane2011a}, we studied a family of Markov processes on
$\partitionsNk$ with the following description.
Let $\nu$ be a finite measure on $\simplexk$. Given an initial state
$\pi\in\partitionsNk$, we construct a Markov process $\boldsymbol
{\Pi
}$ from a
Poisson point process $\mathbf{B}=\{(t,B_t,S_t)\}\subset\mathbb
{R}^+\times\mathcal{P}_{\mathbb{N}:k}^{k}\times{\mathscr{S}_k^k}$ with
intensity $\mathrm{d}t\otimes\varrho_\nu^{\otimes k}\otimes\Upsilon
^{\otimes k}$, where $\Upsilon$ is the uniform distribution on
$\symmetrick$ and, for any measure $\mu$, $\mu^{\otimes k}:=\mu
\otimes\cdots\otimes\mu$ denotes its $k$-fold product measure.
Given a realization of $\mathbf{B}$, we construct $\boldsymbol{\Pi
}:=(\Pi_t, t\geq0)$ from its finite restrictions as follows. First,
for each $n\in\mathbb{N}$, we put $\Pi^{[n]}_0=\pi_{|[n]}$. Then,
for each $t>0$, we write $\Pi^{[n]}_{t-}=(b_1,\ldots,b_r)$, $r\leq
k$, with blocks listed in order of their least element, and
\begin{itemize}
\item if $t>0$ is an atom time of $\mathbf{B}$ with $B_t:=(B^1,\ldots
,B^k)$ a $k$-tuple of partitions and $S_t:=(S_1,\ldots,S_k)$ a
$k$-tuple of permutations $[k]\rightarrow[k]$,
\begin{itemize}
\item we construct the set-valued matrix
%
%e2.6 #&#
%
\begin{equation}
\label{eq:column totals} %
\pmatrix{ B^1_{S_1(1)}\cap
b_1 & B^2_{S_2(1)}\cap b_2 & \cdots&
B^r_{S_r(1)}\cap b_r\vspace*{2pt}
\cr
B^1_{S_1(2)}\cap b_1 & B^2_{S_2(2)}
\cap b_2 & \cdots& B^r_{S_r(2)}\cap b_r
\vspace*{2pt}
\cr
\vdots& \vdots& \ddots& \vdots\vspace*{2pt}
\cr
B^1_{S_1(k)}\cap b_1 & B^2_{S_{2}(k)}
\cap b_2 & \cdots& B^r_{S_r(k)}\cap b_r
} %
,
\end{equation}
and
\item for each $j=1,\ldots, k$, we put $C_j:=\bigcup_{i=1}^r(B^i_{S_i(j)}\cap b_i)$, the union of the entries in row $j$ of
the above matrix. We then define $\Pi^{[n]}_t:=\{C_1,\ldots,C_k\}
\setminus\{\emptyset\}$, provided $\Pi^{[n]}_t\neq\Pi^{[n]}_{t-}$;
\end{itemize}
\item otherwise, we put $\Pi^{[n]}_t=\Pi^{[n]}_{t-}$.
\end{itemize}
We have shown \cite{Crane2011a} that the finite-dimensional transition
rates for this process are
\[
Q_n\bigl(\pi,\pi'\bigr)=k^{\downarrow\#\pi'}\prod
_{b\in\pi}\frac{\varrho
_\nu^{b}(\pi_{|b}')}{k^{\downarrow\#\pi_{|b}'}},\qquad\pi\neq\pi '\in {
\mathcal{P}_{[n]:k}},
\]
where $k^{\downarrow j}:=k(k-1)\cdots(k-j+1)$ and $\varrho_{\nu
}^{b}$ denotes the measure $\varrho_{\nu}$ induces on the space of
partitions of $b\subseteq\mathbb{N}$.

The above construction has an easy description as a three step
procedure. For $k\geq1$, let $\pi:=\{b_1,\ldots,b_r\}$, $r\leq k$,
be a partition of a finite or infinite set. Then we obtain a jump from
$\pi$ to $\pi'$ as follows.\vadjust{\goodbreak}
\begin{enumerate}[(iii)]
\item[(i)] Independently, for each $i=1,\ldots,r$, randomly partition
$b_i$ according to the paintbox process $\varrho_{\nu}$ restricted to
$b_i$. Write $B_i:=\{B_{i,1},\ldots,B_{i,r_i}\}$ to denote the
partition obtained.
\item[(ii)] Independently, for each $i=1,\ldots,r$, randomly label
the blocks of $B_i$ by sampling uniformly without replacement from
$[k]$. Equivalently, we can draw a uniform random permutation $\sigma
_i$ of $[k]$ and order the blocks of $B_i$ by adding $k-r_i$ empty-sets
to the end of $B_i$ and writing $C_i:=(B_{i,\sigma_i(1)},\ldots
,B_{i,\sigma_i(k)})$.
\item[(iii)] We define $\pi'$ by merging all subsets assigned the
same label in step (ii); that is, we put $B_l':=\bigcup_{j=1}^kB_{j,\sigma_j(l)}$ for each $l=1,\ldots,k$ and then define
$\pi':=\{B_1',\ldots,B_k'\}\setminus\{\emptyset\}$.
\end{enumerate}
This procedure produces an exchangeable Feller process on
$\partitionsNk$.
The next example illustrates that an exchangeable Markov process on
$\partitionsNk$ need not be consistent.

%ex2.1 #&#
%
\begin{example}[(Failure of consistency property)]\label{ex:counterexample}
Throughout this example, let $s_0:=(2/3,1/3)\in{\Delta_k^{\downarrow
}}$. With
initial state $\Pi_0\sim\varrho_{s_0}$, we define the infinitesimal
jump rates of $\boldsymbol{\Pi}$ as follows. For every $t\geq0$,
\begin{itemize}
\item given $\Pi_t\neq\mathbf{1}_{\mathbb{N}}$, the trivial
one-block partition of $\mathbb{N}$, $\Pi_t$ jumps to $\mathbf
{1}_{\mathbb{N}}$ at rate $1$, and
\item given $\Pi_t=\mathbf{1}_{\mathbb{N}}$, $\Pi_t$ jumps to
$B\sim\varrho_{s_0}$ at rate 2.
\end{itemize}
Clearly, $\boldsymbol{\Pi}$ is Markovian, exchangeable, and has c\`
adl\`
ag sample paths; however, for each $n\in\mathbb{N}$, the restriction
$\boldsymbol{\Pi}_{|[n]}:=(\Pi_{t|[n]}, t\geq0)$ is not Markovian
because the jump rate at every time $t\geq0$ depends on whether $\Pi
_t$ is trivial, which depends on the tail of $(\Pi_{t|[n]}, n\in
\mathbb{N})$.
\end{example}

We focus on generalizing (i)--(iii). To do so, we work on the space
$\mathL$ of labeled partitions of $\mathbb{N}$ with $k$ classes.
The relationship between $\mathL$ and $\partitionsNk$ is
straightforward, and the added structure of $\mathL$ enables a
cleaner exposition.

%s2.4 #&#
\subsection{Labeled partitions}\label{eq:labeled partitions}

For fixed $k\in\mathbb{N}$, a \emph{$k$-partition} is a labeled set
partition with $k$ classes. Specifically, for any $A\subseteq\mathbb
{N}$, a $k$-partition $\lambda$ of $A$ is a length $k$ set-valued
vector $(\lambda_1,\ldots,\lambda_k)$ with $\lambda_i\subseteq A$
for each $i\in[k]$, $\lambda_i\cap\lambda_{i'}=\emptyset$ for
$i\neq i'$, and $\bigcup_{i=1}^k \lambda_i=A$. Alternatively, for
$A=[n]$, $\lambda$ can be regarded as
\begin{itemize}
\item a sequence $\lambda=\lambda^1\lambda^2\cdots\lambda^n$ in
$[k]^{[n]}$, where
\[
\lambda^i=j \quad\Longleftrightarrow\quad i\in\lambda_j,
\quad \mbox{or}
\]
\item a map $\lambda\dvtx[n]\rightarrow[k]$, where $\lambda
(i)=\lambda
^i$ for each $i\in[n]$.
\end{itemize}
Note that all three specifications of $\lambda$ are equivalent and can
be used interchangeably. In general, we write $\mathcal{L}_{A:k}$ to
denote the space of $k$-partitions of $A\subseteq\mathbb{N}$.

Any $\lambda\in{\mathcal{L}_{[n]:k}} $ induces a
partition of $[n]$ through the map $\mathcal{B}_n\dvtx{\mathcal
{L}_{[n]:k}} \rightarrow{\mathcal
{P}_{[n]:k}} $, defined by
\[
\mathcal{B}_n(\lambda):=\{\lambda_1,\ldots,
\lambda_k\}\setminus\{ \emptyset\},
\]
the unordered collection of classes of $\lambda$ with empty sets
removed. Permutations and injection maps act on $({\mathcal
{L}_{[n]:k}}, n\in\mathbb{N})$ similarly to their action
on $({\mathcal{P}_{[n]:k}}, n\in\mathbb{N})$. In
general, let $\psi\dvtx[m]\rightarrow[n]$, $m\leq n$, be an injection.
Then we define $\psi^*\dvtx{\mathcal{L}_{[n]:k}}
\rightarrow\mathcal{L}_{[m]:k}$ by
\[
\psi^*(\lambda):=\lambda\circ\psi\qquad\mbox{for every }\lambda \in {
\mathcal{L}_{[n]:k}},
\]
where $\lambda\in{\mathcal{L}_{[n]:k}} $ is
treated as a map $[n]\rightarrow[k]$.

The restriction map ${\mathcal{L}_{[n]:k}}
\rightarrow\mathcal{L}_{[m]:k}$ is defined by
\[
\lambda_{|[m]}:=\bigl(\lambda_1\cap[m],\ldots,
\lambda_k\cap[m]\bigr),
\]
and the notion of compatibility for sequences of labeled partitions
carries over from unlabeled partitions.
We define $\mathL$ as the space of $k$-partitions of $\mathbb
{N}$, whose elements can be represented by a compatible sequence of
finite $k$-partitions. Finally, we equip $\mathL$ with ultrametric
%
%e2.7 #&#
%
\begin{equation}
\label{eq:ultrametric}d\bigl(\lambda,\lambda '\bigr):=2^{-\mathbf{n}(\lambda,\lambda')},
\qquad\lambda,\lambda'\in \mathL,
\end{equation}
where $\mathbf{n}(\lambda,\lambda'):=\max\{n\in\mathbb{N}\dvtx
\lambda_{|[n]}=\lambda_{|[n]}'\}$, and $\sigma$-field $\sigma
\langle\bigcup_{n\in\mathbb{N}}{\mathcal
{L}_{[n]:k}}  \rangle$.

The projections $(\mathcal{B}_n, n\in\mathbb{N}), \psi'$, and $\psi
^*$ cooperate with one another; that is, the diagram in \eqref
{eq:commutative diagram} commutes: $\mathcal{B}_m\circ\psi^*=\psi
'\circ\mathcal{B}_n$. By this natural correspondence, we can study
processes on $\mathL$ and later project to $\partitionsNk$. Under
mild conditions, the projection into $\partitionsNk$ preserves most,
and sometimes all, of the properties of a process on $\mathL$.
Using this correspondence, we principally study processes on
$\mathL$ with the intention to later project into $\partitionsNk
$. We discuss this procedure briefly in Section~\ref
{section:partition-processes}, but, by that time, most of its
implications should be obvious.
%
%e2.8 #&#
%
\begin{equation}
\label{eq:commutative diagram} %
\begin{array}{c}
\xymatrix{[n]
\\
[m] \ar[u]^{\psi}} \qquad\quad\xymatrix{\mathcal{L}_{[n]:k}
\ar[d]^{\psi^*} \ar [r]^{\mathcal{B}_n} & \mathcal{P}_{[n]:k}
\ar[d]^{\psi'}
\\
\mathcal{L}_{[m]:k} \ar[r]^{\mathcal{B}_m} & \mathcal{P}_{[m]:k}
} %\begin{tikzpicture}
%\matrix(m) [matrix of math nodes, row sep=3em,
%column sep=3em, text height=1.5ex, text depth=0.25ex]
%{ [n] & & \labelednk&&\partitionsnk\\
% \emm& & \mathcal{L}_{[m]:k} &&\mathcal{P}_{[m]:k} \\};
%\path[->]
%(m-2-1) edge node[auto]{$\psi$} (m-1-1)
%(m-1-3) edge node[auto]{$\psi^*$} (m-2-3)
% edge node[auto]{$\mathcal{B}_n$} (m-1-5)
%(m-2-3) edge node[auto]{$\mathcal{B}_m$} (m-2-5)
%(m-1-5) edge node[auto]{$\psi'$} (m-2-5);
%%edge (m-1-4);
%\end{tikzpicture}
%
\end{array}
\end{equation}

%s2.5 #&#
\subsection{Exchangeable random $k$-partitions}\label{section:random
k-partitions}

A random $k$-partition $\Lambda:=(\Lambda_i, 1\leq i\leq k)$ of
$\mathbb{N}$ is called \emph{exchangeable} if, regarded as a
$[k]$-valued sequence $\Lambda:=\Lambda^1\Lambda^2\cdots$\,, it satisfies
\[
\Lambda^{\sigma}:=\Lambda^{\sigma(1)}\Lambda^{\sigma(2)}\cdots
\equalinlaw\Lambda,
\]
for all permutations $\sigma\dvtx\mathbb{N}\rightarrow\mathbb{N}$
fixing all but finitely many $n\in\mathbb{N}$.

By de Finetti's theorem, the law of an exchangeable $k$-partition is
determined by a unique probability measure $\nu$ on the
$(k-1)$-dimensional simplex
\[
\simplexk:= \Biggl\{(s_1,\ldots,s_k)\dvtx
s_i\geq0\mbox{ and }\sum_{i=1}^ks_i=1
\Biggr\}.
\]
For $s\in\simplexk$, we let $\Lambda^1,\Lambda^2,\ldots$ be i.i.d. from
\[
P_s\bigl\{\Lambda^1=j\bigr\}=s_j,\qquad
j=1,\ldots,k,
\]
and define $\Lambda:=\Lambda^1\Lambda^2\cdots$\,, whose distribution
we denote $\zeta_s$. For a measure $\nu$ on $\simplexk$, we write
\[
\zeta_{\nu}(\cdot):=\int_{\simplexk}\zeta_s(
\cdot)\nu(\mathrm{d}s)
\]
to denote the $\nu$-mixture of $\zeta_s$-measures.

%s3 #&#
\section{Lipschitz partition processes}\label{section:ct Markov}

A random collection $\boldsymbol{\Lambda}=(\Lambda_t, t\geq0)$ in
$\mathL$ is a \emph{Markov process} if, for every $t>0$, the
$\sigma$-fields $\sigma\langle\Lambda_s, s<t\rangle$ and $\sigma
\langle\Lambda_s, s>t\rangle$ are conditionally independent given
$\Lambda_t$.
We are interested in \emph{consistent} Markov processes on $\mathL
$. We specialize to exchangeable processes in Section~\ref
{section:exchangeable processes}.

In Section~\ref{section:coalescent}, we showed a construction of
exchangeable coalescent processes by an iterated application of the
$\Coag$ operator at the atom times of a Poisson point process.
Fundamental properties of the $\Coag$ operator endow the coalescent
process with consistency and the Feller property. Of utmost importance
is Lipschitz continuity, without which the process restricted to, say,
$[n]$ could depend on indices $\{n+1,n+2,\ldots\}$ and the
restrictions need not be Markovian, as in Example~\ref{ex:counterexample}.

%s3.1 #&#
\subsection{Poissonian construction}\label{section:general construction}

Let $\Phi:=\{F\dvtx\mathL\rightarrow\mathL\}$ be the collection
of all maps $\mathL\rightarrow\mathL$ and, for each $n\in
\mathbb{N}$, let $\Phi_n\subseteq\Phi$ be the subcollection of maps
so that the restriction of $F(\lambda)$ to ${\mathcal
{L}_{[n]:k}} $ depends on $\lambda$ only through $\lambda
_{|[n]}$, that is
\[
\Phi_n:=\bigl\{F\in\Phi\dvtx\lambda_{|[n]}=
\lambda_{|[n]}' \Longrightarrow F(\lambda)_{|[n]}=F
\bigl(\lambda'\bigr)_{|[n]}\mbox{ for all }\lambda,\lambda
'\in\mathL\bigr\}.
\]
These collections satisfy
\[
\Phi\supset\cdots\supset\Phi_{n-1}\supset\Phi_n\supset\Phi
_{n+1}\supset\cdots,
\]
whose limit $\bigcap_{n\in\mathbb{N}}\Phi_n=\Phi_{\infty}$ exists
and is non-empty. (For example, the identity map $\idN\dvtx\mathL
\rightarrow\mathL$ is in $\Phi_n$ for every $n\in\mathbb{N}$
and, hence, $\idN\in\Phi_{\infty}$.)
For all $F\in\Phi_{\infty}$, the restriction $F(\lambda)_{|[n]}$
depends only on $\lambda_{|[n]}$, for every $n\in\mathbb{N}$.

%le3.1 #&#
%
\begin{lemma}
The collection $\Phi_\infty$ is in one-to-one correspondence with
\[
\Lip(\mathL):=\bigl\{F\in\Phi\dvtx d\bigl(F(\lambda),F\bigl(\lambda'
\bigr)\bigr)\leq d\bigl(\lambda,\lambda'\bigr)\mbox{ for all }
\lambda,\lambda'\in\mathL \bigr\},
\]
Lipschitz continuous maps $\mathL\rightarrow\mathL$ with
Lipschitz constant $1$.
\end{lemma}

\begin{pf}
First, suppose $F\in\Lip(\mathL)$. Then $d(F(\lambda),F(\lambda
'))\leq d(\lambda,\lambda')$ for every $\lambda,\lambda'\in
\mathL$. By definition of the metric \eqref{eq:ultrametric},
$\lambda_{|[r]}=\lambda_{|[r]}'$ for all $r\leq-\log_2 d(\lambda
,\lambda')$ and $d(F(\lambda),F(\lambda'))\leq2^{-r}$; hence, for
every $n\in\mathbb{N}$, $\lambda_{|[n]}=\lambda_{|[n]}'$ implies
$F(\lambda)_{|[n]}=F(\lambda')_{|[n]}$ and $F\in\Phi_\infty$. The
converse is immediate by the definition of the sets $(\Phi_{n}, n\in
\mathbb{N})$ above.
\end{pf}

As $\Lip(\mathL)$ is exactly the set
\[
\bigl\{F\in\Phi\dvtx\forall n\in\mathbb{N}, \lambda_{|[n]}=\lambda
_{|[n]}' \Longrightarrow F(\lambda)_{|[n]}=F\bigl(
\lambda'\bigr)_{|[n]}\mbox{ for all }\lambda,
\lambda'\in\mathL\bigr\},
\]
any $F\in\Lip(\mathL)$ can be written as the compatible sequence
$(F_{[1]},F_{[2]},\ldots)$ of its restrictions to $\Lip(
{\mathcal{L}_{[n]:k}} )$ for each $n\in\mathbb{N}$.
Specifically, the restriction $F_{[n]}$ of $F\in\Lip(\mathL)$ to
$\Lip({\mathcal{L}_{[n]:k}} )$ is defined, for
every $\lambda\in{\mathcal{L}_{[n]:k}} $, by
$F_{[n]}(\lambda)=F(\lambda^*)_{|[n]}$, for any choice of $\lambda
^*\in\mathL$ such that \mbox{$\lambda^*_{|[n]}=\lambda$}. In this
sense, $\Lip(\mathL)$ is a projective limit space which we can
equip with the ultrametric
%
%e3.1 #&#
%
\begin{equation}
\label{eq:ultrametric lip}{d}_{\Phi
}\bigl(F,F'\bigr):=2^{-\mathbf{n}(F,F')},
\end{equation}
where $\mathbf{n}(F,F'):=\max\{n\in\mathbb{N}\dvtx
F_{[n]}=F_{[n]}'\}$,
and $\sigma$-field $\mathcal{F}=\sigma\langle\bigcup_{n\in\mathbb
{N}}\Lip({\mathcal{L}_{[n]:k}} )\rangle$. It
follows that any measure $\varphi$ on $(\Lip(\mathL),\mathcal
{F})$ determines a measure $\varphi_n$ on $\Lip({\mathcal
{L}_{[n]:k}} )$ through
%
%e3.2 #&#
%
\begin{equation}
\label{eq:restricted measure} \varphi_n(F):=\varphi\bigl(\bigl\{F^*\in\Lip(\mathL)
\dvtx F^*_{[n]}=F\bigr\} \bigr),\qquad F\in\Lip({\mathcal{L}_{[n]:k}}
).
\end{equation}

For any $n\in\mathbb{N}\cup\{\infty\}$, let ${\mathbf
{Id}_n} $ denote the identity map ${\mathcal
{L}_{[n]:k}} \rightarrow{\mathcal
{L}_{[n]:k}} $.\footnote{To maintain consistent notation, we
also define $[\infty]:=\mathbb{N}$ so that ${\mathcal
{L}_{[n]:k}} =\mathL$ for $n=\infty$.} Then, for a
measure $\varphi$ on $(\Lip(\mathL),\mathcal{F})$ satisfying
%
%e3.3 #&#
%
\begin{equation}
\label{eq:finite rates}\varphi\bigl(\{\idN\}\bigr)=0 \quad \mbox{and}\quad
\varphi_n\bigl(\Lip({\mathcal{L}_{[n]:k}} )\setminus\{{
\mathbf{Id}_n} \}\bigr)<\infty\qquad\mbox{for every }n\in\mathbb{N},
\end{equation}
let $\mathbf{F}:=\{(t,F^t)\}\subset\mathbb{R}^+\times\Lip
(\mathL)$ be a Poisson point process with intensity $\mathrm{d}t\otimes
\varphi$. Given $\mathbf{F}$ and some (possibly random) initial state
$\lambda_0\in\mathL$, we construct a Markov process $\boldsymbol
{\Lambda}$ on $\mathL$ as follows. For each $n\in\mathbb{N}$,
we define $\boldsymbol{\Lambda}^{[n]}=(\Lambda^{[n]}_t, t\geq0)$ on
${\mathcal{L}_{[n]:k}} $ by $\Lambda
^{[n]}_0=\lambda_{0|[n]}$ and
%
%e3.4 #&#
%
\begin{equation}
\label{eq:PPP formal} %
\begin{array} {@{}l@{\hspace*{6pt}}l@{}} \bullet&
\mbox{if }t>0\mbox{ is an atom time of }\mathbf{F}\mbox{ such that
}F^t_{[n]}\neq{\mathbf{Id}_n}, \mbox{ we put }
\Lambda^{[n]}_t=F^t_{[n]}\bigl(
\Lambda^{[n]}_{t-}\bigr);
\\
\bullet&\mbox{otherwise, we put }\Lambda^{[n]}_t=
\Lambda^{[n]}_{t-}. \end{array} %
\end{equation}
%
%pr3.1 #&#

\begin{prop}\label{prop:PPP well-def}
For every $n\in\mathbb{N}$, $\boldsymbol{\Lambda}^{[n]}$ is a c\`
adl\`ag
finite state space Markov process, and $(\boldsymbol{\Lambda}^{[n]},
n\in\mathbb{N})$ determines a unique consistent Markov process
$\boldsymbol{\Lambda}$ on $\mathL$.
\end{prop}

\begin{pf}
That each $\boldsymbol{\Lambda}^{[n]}$ is c\`adl\`ag follows from
\eqref
{eq:finite rates} since $\varphi_n(\Lip({\mathcal
{L}_{[n]:k}} )\setminus\{{\mathbf{Id}_n}
\})<\infty$ ensures that, within any bounded interval of $[0,\infty
)$, there are at most finitely many atom times of $\mathbf{F}$ for
which $F^t_{[n]}\neq{\mathbf{Id}_n} $. Furthermore,
for each $n\in\mathbb{N}$, $\boldsymbol{\Lambda}^{[n]}$ is Markov by
the construction in \eqref{eq:PPP formal}. The collection
$(\boldsymbol
{\Lambda}^{[n]}, n\in\mathbb{N})$ is compatible by construction and
therefore, for every $t\geq0$, $({\Lambda}^{[n]}_t, n\in\mathbb
{N})$ determines a unique $k$-partition $\Lambda_t$ of $\mathbb{N}$.
It follows that $(\boldsymbol{\Lambda}^{[n]}, n\in\mathbb{N})$
determines a unique consistent Markov process $\boldsymbol{\Lambda
}=(\Lambda_t, t\geq0)$ on $\mathL$.
\end{pf}

Some remarks about the above construction:
\begin{enumerate}[(iii)]
\item[(i)] $\boldsymbol{\Lambda}$ need not be exchangeable; we treat
exchangeable processes in Section~\ref{section:exchangeable processes}
and give an explicit example of a non-exchangeable process in
Section~\ref{section:non-exchangeable process}.
\item[(ii)] Each restriction $\boldsymbol{\Lambda}_{|[n]}:=(\Lambda
_{t|[n]}, t\geq0)$ has a Poisson point process construction based on
$\mathbf{F}^{(n)}\subset\mathbb{R}^+\times\Lip({\mathcal
{L}_{[n]:k}} )$ with intensity $\mathrm{d}t\otimes\varphi_n$.
\item[(iii)] The second half of \eqref{eq:finite rates} is needed so
that the finite restrictions $(\boldsymbol{\Lambda}_{|[n]}, n\in
\mathbb
{N})$ are c\`adl\`ag. Also, $\varphi$ must put all its support on
$\Lip(\mathL)$, or else the construction in \eqref{eq:PPP
formal} would not result in a compatible collection of finite state
space processes.

The second half of \eqref{eq:finite rates} corresponds to the second
half of \eqref{eq:regularity coalescent} in the following precise
sense. In \eqref{eq:finite rates}, we exclude the identity map
${\mathbf{Id}_n} $ since it does not result in a
jump in the restricted process $\boldsymbol{\Lambda}^{[n]}$, for each
$n\in\mathbb{N}$. Similarly, for each $n\in\mathbb{N}$, $\mathbf
{0}_{[n]}$ is the neutral element for $\Coag$, that is, $\Coag(\pi
,\mathbf{0}_{[n]})=\pi$ for all $\pi\in{\mathcal
{P}_{[n]}} $. Hence, $\mathbf{0}_{[n]}$ determines the
identity map ${\mathcal{P}_{[n]}} \rightarrow
{\mathcal{P}_{[n]}} $ by way of the coagulation operator.
\end{enumerate}

%s3.2 #&#
\subsection{The Feller property}\label{section:Feller property}

Alternatively, we can construct $\boldsymbol{\Lambda}$ from $\mathbf{F}$
by first constructing a Markov process ${\phi}^\infty$ on
$\Lip(\mathL)$. For each $n\in\mathbb{N}$, we construct $\phi
^{[n]}:=(\phi^{[n]}_t, t\geq0)$ on $\Lip({\mathcal
{L}_{[n]:k}} )$ by $\phi^{[n]}_0={\mathbf
{Id}_n} $ and
%
%e3.5 #&#
%
\begin{equation}
\label{eq:PPP formal2} %
\begin{array} {@{}l@{\hspace*{6pt}}l@{}} \bullet&
\mbox{if }t>0\mbox{ is an atom time of }\mathbf{F}\mbox{ such that
}F^t_{[n]}\neq{\mathbf{Id}_n}, \mbox{ we put }
\phi^{[n]}_t=F^t_{[n]}\circ
\phi^{[n]}_{t-};
\\
\bullet&\mbox{otherwise, we put }\phi^{[n]}_t=
\phi^{[n]}_{t-}. \end{array} %
\end{equation}
%
%co3.1 #&#

\begin{cor}\label{cor:lip-ppp}
The collection $(\phi^{[n]}, n\in\mathbb{N})$ is consistent on
$(\Lip({\mathcal{L}_{[n]:k}} ), n\in\mathbb{N})$
and determines a unique Markov process $\phi^\infty$ on $\Lip
(\mathL)$. Moreover,
$\boldsymbol{\Lambda}^{\infty}=(\Lambda_t^{\infty}, t\geq0)$
defined by
%
%e3.6 #&#
%
\begin{equation}
\label{eq:phi rep} \Lambda_t^{\infty}=\phi^\infty_t(
\Lambda_0) \qquad\mbox{for every }t\geq0,
\end{equation}
is a version of $\boldsymbol{\Lambda}$ in \eqref{eq:PPP formal}.
\end{cor}

\begin{pf}
The first claim follows immediately by the arguments in Proposition~\ref{prop:PPP well-def}.

To establish the second claim, let $\mathbf{F}$ be the Poisson point
process with intensity $\mathrm{d}t\otimes\varphi$ and, for every $n\in
\mathbb{N}$, let $J_n$ be the set of atom times of $\mathbf{F}$ such
that $F^t_{[n]}\neq{\mathbf{Id}_n} $. By \eqref
{eq:finite rates}, $J_n\cap[0,t]$ is almost surely finite for every
$n\in\mathbb{N}$ and $t<\infty$. We construct $\boldsymbol{\Lambda}$
from $\mathbf{F}$ as in \eqref{eq:PPP formal} and $\phi^\infty$
from $\mathbf{F}$ as in \eqref{eq:PPP formal2}.

For fixed $n\in\mathbb{N}$ and $t>0$, write $t\geq t_1>\cdots>t_r>0$
to be the ranked atom times of $J_n$ before time $t$.
Each $F^{t_i}\in\Lip(\mathL)$, $i=1,\ldots,r$, and so
\[
\Lambda_{t|[n]}=\bigl(F^{t_1}_{[n]}\circ\cdots\circ
F^{t_r}_{[n]}\bigr) (\Lambda_{0|[n]})=
\phi^\infty_{t,[n]}(\Lambda _{0|[n]})=
\phi^\infty_t(\Lambda_0)_{|[n]}=
\Lambda^\infty_{t|[n]},
\]
where $\phi^{\infty}_{t,[n]}$ denotes the restriction of $\phi
^{\infty}_t$ to $\Lip({\mathcal{L}_{[n]:k}} )$.
Hence, $\boldsymbol{\Lambda}^\infty_{|[n]}=\boldsymbol{\Lambda}_{|[n]}$
almost surely for every $n\in\mathbb{N}$; whence, $\boldsymbol
{\Lambda
}^\infty=\boldsymbol{\Lambda}$ almost surely. The conclusion follows.
\end{pf}
%
%re3.1 #&#

\begin{rmk}
In essence, representation \eqref{eq:phi rep} entails the application
of a flow $(\phi_{s,t}, 0\leq s<t<\infty)$ on the space $\Lip
(\mathL)$, for which we apply $\phi_t:=\phi_{0,t}$ to $\Lambda
_0$, for each $t\geq0$. This can be compared to constructions of
coalescent processes by flows of bridges \cite{BertoinLeGall2005}.
\end{rmk}

Representation \eqref{eq:phi rep} is convenient for studying the
semigroup of $\boldsymbol{\Lambda}$.
For every bounded, continuous function $g\dvtx\mathL\rightarrow
\mathbb{R}$, the semigroup $(\mathbf{P}_t, t\geq0)$ of $\boldsymbol
{\Lambda}$ is defined by
\[
\mathbf{P}_tg(\lambda):=\mathbb{E}_\lambda g(
\Lambda_t),\qquad t\geq0, \lambda\in\mathL,
\]
the expectation of $g(\Lambda_t)$ given $\Lambda_0=\lambda$. In
addition, $(\mathbf{P}_t, t\geq0)$ is called a \emph{Feller
semigroup}, and the process $\boldsymbol{\Lambda}$ is called a \emph
{Feller process}, if, for every bounded, continuous $g\dvtx\mathL
\rightarrow\mathbb{R}$,
\begin{itemize}
\item$\lambda\mapsto\mathbf{P}_t g(\lambda)$ is continuous for
every $t>0$, and
\item$\lim_{t\downarrow0}\mathbf{P}_t g(\lambda)=g(\lambda)$ for
all $\lambda\in\mathL$.
\end{itemize}
%
%co3.2 #&#

\begin{cor}\label{cor:semigroup}
The semigroup $(\mathbf{P}_t, t\geq0)$ of $\boldsymbol{\Lambda}$ satisfies
%
%e3.7 #&#
%
\begin{equation}
\label{eq:semigroup} \mathbf{P}_tg(\lambda):=\mathbb{E}g\bigl(
\phi^\infty_t(\lambda)\bigr),
\end{equation}
for every bounded, continuous map $g\dvtx\mathL\rightarrow
\mathbb
{R}$ and every $\lambda\in\mathL$, where $(\phi^{\infty}_t,
t\geq0)$ is the process in Corollary~\ref{cor:lip-ppp}.
\end{cor}

The proof follows immediately from Corollary~\ref{cor:lip-ppp}.
%
%th3.1 #&#

\begin{them}\label{thm:feller}
The process $\boldsymbol{\Lambda}$ constructed in \eqref{eq:PPP formal}
fulfills the Feller property.
\end{them}

\begin{pf}
Continuity of the map $\lambda\mapsto\mathbf{P}_t g(\lambda)$ is an
immediate consequence of continuity of $g$, the description of $\mathbf
{P}_t$ in \eqref{eq:semigroup}, and the fact that $\phi^\infty_t\in
\Lip(\mathL)$ for all $t>0$ almost surely.

That $\lim_{t\downarrow0}\mathbf{P}_tg(\lambda)=g(\lambda)$ for
all $\lambda\in\mathL$ follows by continuity of $g$ and \eqref
{eq:finite rates}, which ensures that the time of the initial jump out
of $\lambda_{|[n]}$ is strictly positive, for every $n\in\mathbb{N}$.
\end{pf}

By the Feller property, any $\boldsymbol{\Lambda}$ with the construction
in \eqref{eq:PPP formal} has a c\`adl\`ag version.
For the rest of the paper, we implicitly assume $\boldsymbol{\Lambda}$
has c\`adl\`ag paths.
%
%de3.1 #&#

\begin{defn}[(Lipschitz partition process)]\label{def:LPP}
We call the Markov process $\boldsymbol{\Lambda}$ constructed in
\eqref
{eq:PPP formal2} a \emph{Lipschitz partition process} directed by
$\varphi$.
\end{defn}

%s4 #&#
\section{Exchangeable Lipschitz partition processes}\label
{section:exchangeable processes}

A process $\boldsymbol{\Lambda}$ on $\mathL$ is called \emph
{exchangeable} if $\boldsymbol{\Lambda}\equalinlaw\boldsymbol
{\Lambda
}^\sigma$ for all permutations $\sigma\dvtx\mathbb{N}\rightarrow
\mathbb
{N}$ fixing all but finitely many elements of $\mathbb{N}$. We have
already shown (Proposition~\ref{prop:PPP well-def} and Theorem~\ref
{thm:feller}) that Lipschitz partition processes are consistent and
possess the Feller property.
We now consider \emph{exchangeable Lipschitz partition processes} on
$\mathL$.

Provided its rate measure $\mu$ is exchangeable, a coalescent process
(Section~\ref{section:coalescent}) is exchangeable. In the
exchangeable case, the directing measure $\mu$ in \eqref
{eq:regularity coalescent} need only satisfy $\mu(1\sim2)<\infty$.
Furthermore, if we describe $\mu$ by a paintbox measure $\varrho_{\nu
}$ on $\partitionsN$, \eqref{eq:regularity coalescent} implies
\[
\nu\bigl(\bigl\{(0,0,\ldots)\bigr\}\bigr)=0 \quad\mbox{and}\quad\int
_{{\Delta
^{\downarrow}}}(1-s_1)\nu (\mathrm{d}s)<\infty.
\]
For Lipschitz partition processes constructed in \eqref{eq:PPP
formal}, $\varphi$ must be restricted to the space of \emph{strongly
Lipschitz maps} to ensure exchangeability. We introduce strongly
Lipschitz maps in Section~\ref{section:strong lipschitz} and show some
of their properties in Section~\ref{section:matrix multiplication}.

%s4.1 #&#
\subsection{Strongly Lipschitz maps}\label{section:strong lipschitz}

In this section, we see that any exchangeable Markov process
$\boldsymbol
{\Lambda}$ with construction \eqref{eq:PPP formal} must be directed
by a measure $\varphi$ whose support is contained in the proper subset
of \emph{strongly Lipschitz maps} on $\mathL$.

For any $A\subseteq\mathbb{N}$ and $\lambda,\lambda'\in\mathcal
{L}_{A:k}$, we define the \emph{overlap} of $\lambda$ and $\lambda'$ by
%
%e4.1 #&#
%
\begin{equation}
\label{eq:inner product}\lambda\cap\lambda ':=\bigcup
_{i=1}^k\bigl(\lambda_i\cap
\lambda_i'\bigr),
\end{equation}
and let
%
%e4.2 #&#
%
\begin{equation}
\label{eq:superlipschitz maps n}\Sigma_n:=\bigl\{F\in\Lip (\mathL)\dvtx
F_{[n]}(\lambda)\cap F_{[n]}\bigl(\lambda'\bigr)
\supseteq \lambda \cap\lambda'\mbox{ for all }\lambda,
\lambda'\in{\mathcal {L}_{[n]:k}} \bigr\}
\end{equation}
be the subset of functions $F\in\Lip(\mathL)$ for which the
overlap of the image of any $\lambda,\lambda'\in{\mathcal
{L}_{[n]:k}} $ by the restriction $F_{[n]}$ contains the
overlap of $\lambda$ and $\lambda'$. By definition of the ultrametric
\eqref{eq:ultrametric} on $\mathL$, if $d(\lambda,\lambda')\leq
2^{-n}$ for some $n\in\mathbb{N}$, then $[n]\subseteq\lambda\cap
\lambda'$; thus, $\Sigma_n\subseteq\Phi_n$ for all $n\in\mathbb
{N}$. We write $\Sigma:=\bigcap_{n\in\mathbb{N}}\Sigma_n$ to
denote the collection of Lipschitz continuous maps satisfying
%
%e4.3 #&#
%
\begin{equation}
\label{eq:superlipschitz maps}F(\lambda)\cap F\bigl(\lambda'\bigr)\supseteq\lambda
\cap\lambda'\qquad\mbox{for all }\lambda ,\lambda'\in
\mathL,
\end{equation}
and we call any $F\in\Sigma$ \emph{strongly Lipschitz continuous}.
In the following proposition, let $\boldsymbol{\Lambda}$ be a Lipschitz
partition process directed by $\varphi$.
%
%pr4.1 #&#

\begin{prop}\label{thm:exchangeable lipschitz}
If $\boldsymbol{\Lambda}$ is exchangeable, then $\varphi$ is supported
on $\mathcal{F}\cap\Sigma$, the trace $\sigma$-field of $\mathcal
{F}=\sigma\langle\bigcup_{n\in\mathbb{N}}\Lip({\mathcal
{L}_{[n]:k}} )\rangle$.
\end{prop}

\begin{pf}
Suppose $\boldsymbol{\Lambda}$ is exchangeable and fix $n\in\mathbb
{N}$. Then $\boldsymbol{\Lambda}^\sigma_{|[n]}\equalinlaw
\boldsymbol
{\Lambda}_{|[n]}$ for all $\sigma\in\symmetricn$. Hence, we can
construct $\boldsymbol{\Lambda}^\sigma$ and $\boldsymbol{\Lambda}$ from
the same Poisson point process $\mathbf{F}:=\{(t,F^t)\}$ by putting
%
%e4.4 #&#
%
\begin{equation}
\label{eq:simultaneous construction}\Lambda^\sigma _{t|[n]}=\sigma^*F^t_{[n]}(
\Lambda_{t-|[n]})=\sigma^*F^t_{[n]}\sigma
^{*^{-1}}\bigl(\Lambda^\sigma_{t-|[n]}\bigr)
\end{equation}
for every $t\in J_n:=\{t>0\dvtx (t,F^t)\in\mathbf{F}\mbox{ and
}F^t_{[n]}\neq{\mathbf{Id}_n} \}$, the jump times
of $\boldsymbol{\Lambda}_{|[n]}$. By \eqref{eq:simultaneous
construction} and the construction of $\boldsymbol{\Lambda}$ in
\eqref
{eq:PPP formal}, $\mathbf{F}^\sigma:=\{(t,\sigma^*F^t\sigma
^{*^{-1}})\}$ has the same law as a Poisson point process on $[0,\infty
)\times\Lip(\mathL)$ with intensity $\mathrm{d}t\otimes\varphi$, for
all $\sigma\in\symmetricn$. Since we have assumed that the support
of $\varphi$ is a subset of $\Lip(\mathL)$ and the set $J_\infty
$ of atom times of $\mathbf{F}$ is at most countable by \eqref
{eq:finite rates}, we have $\sigma^*F^t\sigma^{*^{-1}}\in\Lip
(\mathL)$ for all $t\in J_\infty$ almost surely. It follows that
$\varphi$ must be supported on
\[
\Xi:=\bigl\{F\in\Lip(\mathL)\dvtx\sigma^*F\sigma^{*^{-1}}\in\Lip (\mathL)
\mbox{ for every finite permutation }\sigma\dvtx\mathbb {N}\rightarrow\mathbb{N}
\bigr\},
\]
which is non-empty. To see that $\Xi\subset\Sigma$, we need the
following lemma.
%
%le4.1 #&#

\begin{lemma}\label{lemma:perm}
For $n\in\mathbb{N}$, let $\lambda,\lambda'\in{\mathcal
{L}_{[n]:k}} $ have overlap of size $\#(\lambda\cap\lambda
')=r\in[n]$. Then there exists $\sigma\in\symmetricn$ such that
$\sigma^2$ is the identity $[n]\rightarrow[n]$ and $\lambda^\sigma
\cap\lambda^{\prime\sigma}=[r]$.
\end{lemma}

\begin{pf}
For $m,m'\leq r$, let $r<i_1<i_2<\cdots<i_m$ be the elements of
$(\lambda\cap\lambda')\setminus[r]$ and let $j_1<\cdots<j_{m'}\leq
r$ be the elements of $(\lambda\cap\lambda')^{c}\cap[r]$. Note that
$m'=r-(r-m)=m$; so we can define $\sigma\in\symmetricn$ by $\sigma
(i_l)=j_l$ and $\sigma(j_l)=i_l$ for every $l=1,\ldots,m$, and
$\sigma(i)=i$ otherwise. Clearly, $\sigma^2$ is the identity and
$i\in\lambda\cap\lambda'$ implies $\sigma(i)\in[r]$.
\end{pf}

Now, fix $n\in\mathbb{N}$ and take $F\in\Xi$. For any $\sigma\in
\symmetricn$, we write $F^{\sigma}_{[n]}:=\sigma^{*}F_{[n]}\sigma
^{*^{-1}}$. Take any $\lambda,\lambda'\in{\mathcal
{L}_{[n]:k}} $ and let $\sigma$ be the permutation of $[n]$
from the preceding lemma. Then $\sigma^*=\sigma^{*^{-1}}$, $F^\sigma
_{[n]}:=\sigma^*F_{[n]}\sigma^*\in\Lip({\mathcal
{L}_{[n]:k}} )$, and $F_{[n]}=\sigma^*F^\sigma_{[n]}\sigma
^*$. Let $d_n$ denote the restriction of the metric $d$ in \eqref
{eq:ultrametric} to ${\mathcal{L}_{[n]:k}} $. By
Lipschitz continuity and Lemma~\ref{lemma:perm},
\[
d_n\bigl(F^\sigma_{[n]}\sigma^*(
\lambda),F^\sigma_{[n]}\sigma^*\bigl(\lambda '
\bigr)\bigr)=d_n\bigl(F^\sigma_{[n]}\bigl(
\lambda^\sigma\bigr),F^\sigma_{[n]}\bigl(\lambda
^{\prime\sigma}\bigr)\bigr)\leq2^{-r};
\]
hence, $\lambda^\sigma(j)=\lambda^{\prime\sigma}(j)$ and $[F^\sigma
_{[n]}\sigma^*(\lambda)](j)=[F^\sigma_{[n]}\sigma^*(\lambda')](j)$
for all $j\in[r]$. Finally, take $i\in\lambda\cap\lambda'$. Then
$\sigma(i)\in[r]$ by Lemma~\ref{lemma:perm}, which implies
%
%\[
%
\begin{eqnarray*}
\bigl[F_{[n]}(\lambda)\bigr](i)&=&\sigma^*\bigl[F^\sigma_{[n]}
\sigma^*(\lambda )\bigr](i)
\\
&=&\bigl[F^\sigma_{[n]}\sigma^*(\lambda)\bigr]\bigl(\sigma(i)
\bigr)
\\
&=&\bigl[F^\sigma_{[n]}\sigma^*\bigl(\lambda'
\bigr)\bigr]\bigl(\sigma(i)\bigr)
\\
&=&\sigma^*\bigl[F^\sigma_{[n]}\sigma^*\bigl(
\lambda'\bigr)\bigr](i)=\bigl[F_{[n]}\bigl(
\lambda'\bigr)\bigr](i),
\end{eqnarray*}
%
%\]
%
and $i\in F(\lambda)\cap F(\lambda')$. It follows that $\Xi\subset
\Sigma$.
\end{pf}
%
%re4.1 #&#

\begin{rmk}
The converse of Proposition~\ref{thm:exchangeable
lipschitz} does not hold.
\end{rmk}

Proposition~\ref{thm:exchangeable lipschitz} shows that the directing
measure of an exchangeable Lipschitz partition process can only assign
positive measure to events in the trace $\sigma$-field $\mathcal
{F}\cap\Sigma$.
In the next section, we use condition \eqref{eq:superlipschitz maps}
to characterize the space $\Sigma$.

%s4.2 #&#
\subsection{Strongly Lipschitz maps and set-valued matrix
multiplication}\label{section:matrix multiplication}

A $k\times k$ matrix $M$ over $S\subset\mathbb{N}$ is a collection
$(M_{ij}, 1\leq i,j\leq k)$ of subsets of $S$ for which we define the
operation \emph{multiplication} by
%
%e4.5 #&#
%
\begin{equation}
\label{eq:matrix multiplication}\bigl(M*M'\bigr)_{ij}\equiv
\bigl(MM'\bigr)_{ij}:=\bigcup
_{l=1}^k\bigl(M_{il}\cap
M_{lj}'\bigr),\qquad1\leq i,j\leq k.
\end{equation}
The operation in \eqref{eq:matrix multiplication} mimics
multiplication of real-valued matrices, but for matrices taking values
in a distributive lattice. Here, the lattice operations $\cap$ and
$\cup$ correspond to multiplication and addition, respectively.

We are particularly interested in \emph{partition operators}, matrices
$M$ over $[n]$ with each $M^j\in\mathcal{L}_{[n]:k}$, $j=1,\ldots
,k$, where $M^j$ denotes the $j$th column of $M$.
We write ${\mathcal{M}_{[n]:k}} $ to denote the set
of $k\times k$ partition operators over $[n]$.

Some observations about the operation \eqref{eq:matrix multiplication}:
\begin{enumerate}[(iii)]
\item[(i)] For $m\leq n$, we can define the restriction of $M\in
{\mathcal{M}_{[n]:k}} $ to $\mathcal{M}_{[m]:k}$.
First, we let $I_m^k:=\diag([m],\ldots,[m])$ be the $k\times k$
matrix with diagonal entries $[m]$ and off-diagonal entries the empty set.
Then, for any $M\in{\mathcal{M}_{[n]:k}} $, the
product $M_{[m]}:=I_m^kM=MI_m^k\in\mathcal{M}_{[m]:k}$ is
well-defined as the restriction of $M$ to $\mathcal{M}_{[m]:k}$. It follows
that $({\mathcal{M}_{[n]:k}}, n\in\mathbb{N})$ is
a projective system with limit space $\matrixNk$, partition operators
on $\mathL$.
\item[(ii)] For any injection $\psi:=(\psi_1,\ldots,\psi
_k)\dvtx[m]^k\rightarrow[n]^k$, $m\leq n$, we define the projection
$\psi
^{**}\dvtx{\mathcal{M}_{[n]:k}} \rightarrow\mathcal
{M}_{[m]:k}$ by
\[
\psi^{**}(M):=\bigl(\psi^*_1\bigl(M^1\bigr),
\ldots,\psi^*_k\bigl(M^k\bigr)\bigr), \qquad\mbox{for
every }M\in{\mathcal{M}_{[n]:k}},
\]
where we write $M:=(M^1,\ldots,M^k)$ as the vector of its columns.
In particular, for $\sigma\in\symmetricn$, we write $\sigma
^{**}M=M^\sigma=(\sigma^*M^1,\ldots,\sigma^*M^k)$, the image of $M$
under relabeling each of its columns by $\sigma$.
\item[(iii)] We can equip $\matrixNk$ with the ultrametric ${d}_{\Phi
}$ in \eqref{eq:ultrametric lip} restricted to $\matrixNk$; in particular,
\[
d_{\Phi}\bigl(M,M'\bigr):=2^{-\mathbf{n}(M,M')},
\]
where $\mathbf{n}(M,M'):=\max\{n\in\mathbb{N}\dvtx MI_n^k=M'I_n^k\}$.
\end{enumerate}
We record some facts about partition operators.
%
%le4.2 #&#

\begin{lemma}\label{prop:routine}
Let $n\in\mathbb{N}\cup\{\infty\}$ and $m\leq n$.
\begin{enumerate}[(iii)]
\item[(i)] Any $M\in\mathcal{M}_{[m]:k}$ determines a map $M\dvtx
{\mathcal{M}_{[n]:k}} \rightarrow\mathcal{M}_{[m]:k}$,
$M'\mapsto MM'$.
\item[(ii)] Any $M\in\mathcal{M}_{[m]:k}$ determines a map
$M\dvtx{\mathcal{L}_{[n]:k}} \rightarrow\mathcal
{L}_{[m]:k}$ by
%
%e4.6 #&#
%
\begin{equation}
\label{eq:matrix map}(M\lambda)_i:=\bigcup_{j=1}^k(M_{ij}
\cap\lambda_j),\qquad i=1,\ldots,k, \lambda\in {\mathcal{L}_{[n]:k}}.
\end{equation}
\item[(iii)] The operation \eqref{eq:matrix multiplication} is
associative, that is, $M(M'M'')=(MM')M''$ for all $M,M',M''\in
{\mathcal{M}_{[n]:k}} $.
\item[(iv)] Each $M\in{\mathcal{M}_{[n]:k}} $
determines a Lipschitz continuous map $M\dvtx{\mathcal
{M}_{[n]:k}} \rightarrow{\mathcal
{M}_{[n]:k}} $ through \eqref{eq:matrix multiplication} and
$M\dvtx{\mathcal{L}_{[n]:k}} \rightarrow
{\mathcal{L}_{[n]:k}} $ through \eqref{eq:matrix map}.
\end{enumerate}
\end{lemma}

\begin{pf}
The proof is routine, but we include the proof of (iv) because it is
crucial to the paper. Note that the restriction of any $\lambda\in
\mathL$ to $n\in\mathbb{N}$ can be expressed as $\lambda
_{|[n]}=I_n^k\lambda$. Let $\lambda,\lambda'\in\mathL$ be such
that $ I_r^k\lambda=I_r^k\lambda'$ for some $r\in\mathbb{N}$. Then
$d(\lambda,\lambda')\leq2^{-r}$ and, for every $M\in\matrixNk$,
\[
I_r^k(M\lambda)=\bigl(I_r^kM
\bigr)\lambda=\bigl(MI_r^k\bigr)\lambda=M
\bigl(I_r^k\lambda \bigr)=M\bigl(I_r^k
\lambda'\bigr)=I_r^k\bigl(M
\lambda'\bigr),
\]
implying $d(M\lambda,M\lambda')\leq d(\lambda,\lambda')$.
\end{pf}

%ex4.1 #&#
%
\begin{example}[(Partition operator)]\label{ex:partition operator}
Fix $n=6$, $k=2$, and let $\lambda=(\{1,3,4,5\},\{2,6\})$. Then the
image of $\lambda$ by
\[
M:= %
\pmatrix{ \{2,3\} & \{2,4,5,6\}
\cr
\{1,4,5,6\} & \{1,3\} }
\]
is
\begin{eqnarray*}
M\lambda&:=& %
\pmatrix{ \{2,3\} & \{2,4,5,6\}\vspace*{2pt}
\cr
\{1,4,5,6\} & \{1,3
\} } %
\pmatrix{ \{1,3,4,5\}\vspace*{2pt}
\cr
\{2,6\} } %
\\
&=& %
\pmatrix{ \bigl(\{2,3\}\cap\{1,3,4,5\}\bigr)\cup\bigl(\{2,4,5,6
\}\cap\{2,6\}\bigr)\vspace*{2pt}
\cr
\bigl(\{1,4,5,6\}\cap\{1,3,4,5\}\bigr)\cup\bigl(\{1,3\}\cap
\{2,6\}\bigr) } %
\\
&=& %
\pmatrix{ \{2,3,6\}\vspace*{2pt}
\cr
\{1,4,5\} } %
.
\end{eqnarray*}
\end{example}

%re4.2 #&#

\begin{rmk}[(Partition operators and the $\operatorname{\mathbf{Coag}}$ operator)]
There is a relationship between partition operators and the coagulation
operator $\Coag\dvtx\partitionsN\times\partitionsN\rightarrow
\partitionsN$ from Section~\ref{section:coalescent}. For $k\in
\mathbb{N}$, let $\pi:=\{b_1,\ldots,b_k\}\in\partitionsNk$ and
define $\lambda:=(b_1,\ldots,b_k)$, the $k$-partition obtained by
listing the blocks of $\pi$ in ascending order of their least element.
Now, given $\pi'=\{b_1',\ldots,b_{r'}'\}\in\mathcal{P}_{[k]}$, we
define $M:=M_{\pi'}$ by
\[
M_{ij}:=\cases{ %
%\begin{array}{cc}
\mathbb{N},&\quad$j\in
b_i'$,
\cr
\emptyset,&\quad$\mbox{otherwise}$.
%\end{array}
%%
%\right.
}
\]
Then $\mathcal{B}_{\infty}(M_{\pi'}\lambda)=\Coag(\pi,\pi')$.
For example, let $\pi=123/45/678/9$ so that $\lambda=(123,45,\allowbreak  678,9)$,
and let $\pi'=12/34$.
In this case, $\Coag(\pi,\pi')=12345/6789$ and
\[
M_{\pi'}\lambda= %
\pmatrix{ \mathbb{N} & \mathbb{N} &
\emptyset& \emptyset
\cr
\emptyset& \emptyset& \mathbb{N} & \mathbb{N}
\cr
\emptyset& \emptyset& \emptyset& \emptyset
\cr
\emptyset& \emptyset& \emptyset&
\emptyset } %
\pmatrix{ 123
\cr
45
\cr
678
\cr
9 } %
= %
\pmatrix{ 12345
\cr
6789
\cr
\emptyset
\cr
\emptyset } %
.
\]
Note that, in general, partition operators cannot be used instead of
the coagulation operator in the construction of the coalescent process
because, in the standard coalescent, the initial state $\Pi_0:=\mathbf
{0}_{\mathbb{N}}$ has infinitely many blocks, but partition operators
are defined as $k\times k$ matrices for finite $k\geq1$.
\end{rmk}

%pr4.2 #&#

\begin{prop}\label{cor:isomorphism}
The space $\matrixNk$ of partition operators is in one-to-one
correspondence with $\Sigma$ defined in \eqref{eq:superlipschitz maps n}.
\end{prop}

\begin{pf}
Let $F\in\Sigma$ and $n\in\mathbb{N}$. Then $F\in\Sigma_n$ and,
for each $i\in[n]$, if $\lambda(i)=\lambda'(i)$ then
$F_{[n]}(\lambda_{|[n]})(i)=F_{[n]}(\lambda_{|[n]}')(i)$. For
$j=1,\ldots,k$, let $E_j^{(n)}\in{\mathcal
{L}_{[n]:k}} $ be the $k$-partition of $[n]$ satisfying
$E_j^{(n)}(i)=j$ for every $i\in[n]$. Construct $M_{[n]}\in
{\mathcal{M}_{[n]:k}} $ by setting its $j$th column
$M^j_{[n]}$ equal to the image of $E_j^{(n)}$ by $F_{[n]}$. So
$M_{[n]}:=(F_{[n]}(E_1^{(n)}),F_{[n]}(E_2^{(n)}),\ldots
,F_{[n]}(E_k^{(n)}))$. By definition of $\Sigma_n$ in \eqref
{eq:superlipschitz maps n}, it is clear that $M_{[n]}\lambda
=F_{[n]}(\lambda)$, for every $\lambda\in{\mathcal
{L}_{[n]:k}} $. The collection $(M_{[n]}, n\in\mathbb{N})$
is compatible with respect to the restriction maps on $(
{\mathcal{M}_{[n]:k}}, n\in\mathbb{N})$ and therefore
determines a unique $M\in\matrixNk$ satisfying
\[
M\lambda=F(\lambda) \qquad\mbox{for every }\lambda\in\mathL.
\]

The opposite morphism $\matrixNk\rightarrow\Sigma$ follows from
definition \eqref{eq:matrix map} and definition of the metric in
\eqref{eq:ultrametric}.
\end{pf}

From Proposition~\ref{cor:isomorphism}, we can assume, without loss of
generality, that any exchangeable process with construction \eqref
{eq:PPP formal} is directed by a measure $\mu$ on $ (\matrixNk
,\sigma \langle\bigcup_{n\in\mathbb{N}}{\mathcal
{M}_{[n]:k}}  \rangle  )$ for which
%
%e4.7 #&#
%
\begin{equation}
\label{eq:finite rates matrix}\mu\bigl(\bigl\{I_\infty^k\bigr\}\bigr)=0 \quad
\mbox{and}\quad\mu_n\bigl({\mathcal{M}_{[n]:k}} \setminus \bigl
\{I_n^k\bigr\}\bigr)<\infty\qquad\mbox{for all }n\in
\mathbb{N},
\end{equation}
where $I_{n}^k$ is the partition operator with diagonal entries $[n]$
and off-diagonal entries the empty set, and $\mu_n$ denotes the
restriction of $\mu$ to ${\mathcal{M}_{[n]:k}} $.
Note that \eqref{eq:finite rates matrix} agrees with \eqref{eq:finite rates}.

%th4.1 #&#

\begin{them}\label{cor:exchangeable lipschitz}
Let $\boldsymbol{\Lambda}:=(\Lambda_t, t\geq0)$ be a Lipschitz
partition process on $\mathL$. Then $\boldsymbol{\Lambda}$ is
exchangeable if and only if its directing measure $\mu$
\begin{itemize}
\item is supported on $\matrixNk$,
\item satisfies
%
%e4.8 #&#
%
\begin{equation}
\label{eq:exch finite rates matrix} \mu\bigl(\bigl\{I_{\infty}^k\bigr\}\bigr)=0
\quad\mbox{and}\quad\mu_2\bigl(\bigl\{M\in\mathcal {M}_{[2]:k}
\dvtx M\neq I_2^k\bigr\}\bigr)<\infty,
\end{equation}
and
\item for every permutation $\sigma\dvtx\mathbb{N}\rightarrow
\mathbb
{N}$ fixing all but finitely many $n\in\mathbb{N}$ and every
measurable subset $A\subseteq\matrixNk$,
%
%e4.9 #&#
%
\begin{equation}
\label{eq:30}\mu(A)=\mu\bigl(\bigl\{M^\sigma\dvtx M\in A\bigr\}\bigr).
\end{equation}
\end{itemize}
\end{them}
\begin{pf}
Support of $\mu$ on $\matrixNk$ follows from Proposition~\ref
{cor:isomorphism}, and \eqref{eq:30} is a consequence of
exchangeability and the fact that, for any $M\in{\mathcal
{M}_{[n]:k}} $, $\lambda\in{\mathcal
{L}_{[n]:k}} $, and $\sigma\in\symmetricn$, $(M\lambda
)^\sigma=M^\sigma\lambda^\sigma$. Condition \eqref{eq:exch finite
rates matrix} follows from \eqref{eq:finite rates matrix}.

To show the converse, we need only show that \eqref{eq:exch finite
rates matrix} implies \eqref{eq:finite rates matrix}.
Indeed, let $n\in\mathbb{N}$ and note that the event $
{\mathcal{M}_{[n]:k}} \setminus\{I_n^k\}=\{M\in
{\mathcal{M}_{[n]:k}} \dvtx M\neq I_n^k\}$ implies that there is
some permutation $\sigma\in\symmetricn$ such that $M^{\sigma
}I_2^k\neq I_2^k$; hence,
\begin{eqnarray*}
\mu_n\bigl(\bigl\{M\dvtx M\neq I_n^k\bigr\}
\bigr)&=&\mu_n \biggl(\bigcup_{\sigma\in
\symmetricn} \bigl
\{M\dvtx M^{\sigma}I_2^k\neq I_2^k
\bigr\} \biggr)
\\
&\leq&\sum_{\sigma\in\symmetricn}\mu_2\bigl(\bigl\{M\in
\mathcal {M}_{[2]:k}\dvtx M\neq I_2^k\bigr\}\bigr)
\\
&=&n!\mu_2\bigl(\bigl\{M\in\mathcal{M}_{[2]:k}\dvtx M\neq
I_2^k\bigr\}\bigr)
\\
&<&\infty.
\end{eqnarray*}
The rest is immediate.
\end{pf}

Proposition~\ref{cor:isomorphism} and Theorem~\ref{cor:exchangeable
lipschitz} suggest a construction of arbitrary partition operators. In
brief, take any collection $\lambda^{(1)},\ldots,\lambda^{(k)}$ in
$\mathL$ and, for each $j=1,\ldots,k$, put the $j$th column of
$M\in\matrixNk$ equal to $\lambda^{(j)}$. Likewise, a measure $\mu$
on $\matrixNk$ can be defined by a measure on the product space
${\mathcal{L}^k_{\mathbb{N}:k}}$. Furthermore, using the above
observation, we can
construct a measure $\varphi$ with support in $\Lip(\mathL)$ but
not in $\matrixNk$, leading to an explicit construction of a
non-exchangeable Lipschitz process whose semigroup is not determined by
strongly Lipschitz functions. We show such a process in Section~\ref
{section:non-exchangeable process}.

%s4.3 #&#
\subsection{Examples: Exchangeable Lipschitz partition
processes}\label{section:examples}

%ex4.2 #&#
%
\begin{example}[(Self-similar exchangeable Markov process on
$\boldsymbol{\mathcal{L}}_{\boldsymbol{\mathbb{N}:k}}$)]\label{ex:self-similar2}
For any probability measure $\nu$ on $\simplexk$, recall the
definition of $\zeta_{\nu}$ (Section~\ref{section:random k-partitions}).
Given a measure $\nu$ on $\simplexk$, we write $\mu_{\nu^{\otimes
k}}$ to denote the measure on $\matrixNk$ coinciding with the product
measure $\zeta_\nu\otimes\cdots\otimes\zeta_{\nu}$ on
${\mathcal{L}^k_{\mathbb{N}:k}}$. More generally, for measures $\nu
_1,\ldots,\nu_k$ on
$\simplexk$, $\mu_{\nu_1\otimes\cdots\otimes\nu_k}$ is the
measure on $\matrixNk$ coinciding with $\zeta_{\nu_1}\otimes\cdots
\otimes\zeta_{\nu_k}$ on ${\mathcal{L}^k_{\mathbb{N}:k}}$.

Let $\nu_1,\ldots,\nu_k$ be measures on $\simplexk$ such that
\[
\int_{\simplexk}(1-s_i)\nu_i(\mathrm{d}s)<
\infty\qquad\mbox{for all } i=1,\ldots,k.
\]
Then the second half of \eqref{eq:finite rates} is satisfied for $\mu
_{\nu_1\otimes\cdots\otimes\nu_k}$ and we can construct a process
$\boldsymbol{\Lambda}:=(\Lambda_t, t\geq0)$ from a Poisson point
process $\mathbf{M}:=\{(t,M_t)\}\subset\mathbb{R}^+\times\matrixNk
$ with intensity $\mathrm{d}t\otimes\mu_{\nu_1\otimes\cdots\otimes\nu
_k}$, just as in \eqref{eq:PPP formal}. The infinitesimal jump rates
of this process are given explicitly by
\[
Q_n\bigl(\lambda,\lambda'\bigr):=\prod
_{i=1}^k\zeta_{\nu_i}^{\lambda
_i}\bigl(
\lambda_{|\lambda_i}'\bigr),\qquad\lambda\neq
\lambda'\in {\mathcal{L}_{[n]:k}},
\]
for each $n\in\mathbb{N}$, where $\zeta_{\nu}^b$ denotes the
measure induced on $\mathcal{L}_{b:k}$ by $\zeta_{\nu}$ for any
$b\subseteq\mathbb{N}$. This process is the analog of the
self-similar processes in Section~\ref{section:self-similar}.
\end{example}

%ex4.3 #&#
%
\begin{example}\label{section:one column}
Similar to the above example, let $\nu$ be a measure on $\simplexk$
so that
\[
\zeta_\nu^{(n)}\bigl({\mathcal{L}_{[n]:k}}
\backslash\bigl\{ E_i^{(n)}\bigr\}\bigr)<\infty, \qquad
\mbox{for every }n\in\mathbb{N}\mbox{ and all }i=1,\ldots,k,
\]
where $E_i^{(n)}\in{\mathcal{L}_{[n]:k}} $ is the
$k$-partition of $[n]$ with all elements labeled $i$. With $\mathscr
{U}_k$ denoting the uniform distribution on $[k]$, the Poisson point
process $\mathbf{F}=\{(t,\lambda_t,U_t)\}\subset[0,\infty)\times
\mathL\times[k]$, with intensity $\mathrm{d}t\otimes\zeta_\nu\otimes
\mathscr{U}_k$, determines a random subset $\mathbf{M}\subset
[0,\infty)\times\matrixNk$, where for each atom time $t>0$ of
$\mathbf{F}$ we define $M_t\in\matrixNk$ by putting
\[
M_t^i=\cases{ %
%\begin{array}{cc}
\lambda_t, &\quad$i=U_t$,
\cr
E_i, & \quad$
\mbox{otherwise}$; %\end{array}
%%
%\right.
}
\]
that is, writing $\lambda_t=(\lambda_{t,1},\ldots,\lambda_{t,k})$,
we put
\[
M_t:=\bordermatrix{ & 1 & 2 & \cdots& U_t & \cdots& k
\cr
& \mathbb{N} & \emptyset& \cdots& \lambda_{t,1} & \cdots& \emptyset
\cr
& \emptyset& \mathbb{N} & \cdots& \lambda_{t,2} & \cdots&\emptyset
\cr
& \vdots& \vdots& \ddots& \vdots& \ddots& \vdots
\cr
& \emptyset& \emptyset&
\cdots& \lambda_{t,k} & \cdots& \mathbb{N}}.
\]
Given $\mathbf{F}$ and an initial state $\Lambda_0\in\mathL$,
we construct the process $\boldsymbol{\Lambda}$ as in \eqref{eq:PPP
formal} by putting $\Lambda_t=M_t\Lambda_{t-}$ whenever $t>0$ is an
atom time of $\mathbf{F}$. Variations of this description, for
example, for
which at most one class of the current state $\Lambda_t$ is broken
apart in any single jump, are possible and straightforward. For
example, the rates at which different classes experience jumps need not
be identical.
\end{example}

%ex4.4 #&#
%
\begin{example}[(Group action on $\boldsymbol{\mathcal{L}}_{\boldsymbol{\mathbb{N}:k}}$)]\label{subsec:group action}
For any $\lambda\in\mathL$, we define $\bM_\lambda\in
\matrixNk$ by
%
%e4.10 #&#
%
\begin{equation}
\label{eq:partition to matrix} \bM_\lambda:= %
\pmatrix{ \lambda_1
& \lambda_k & \lambda_{k-1} & \cdots& \lambda_2
\cr
\lambda_2 & \lambda_1 & \lambda_k &
\cdots& \lambda_3
\cr
\vdots& \vdots& \vdots& \ddots& \vdots
\cr
\lambda_k & \lambda_{k-1} & \lambda_{k-2} &
\cdots& \lambda_1 } %
.
\end{equation}
In words, $\bM_\lambda$ is the $k\times k$ matrix whose $j$th column
is the $j$th cyclic shift of the classes of $\lambda$. Note that $\bM
_\lambda\lambda'=\bM_{\lambda'} \lambda$ for all $\lambda,\lambda
'\in\mathL$.

Any measure $\zeta$ on $\mathL$ determines a measure $\mu_{\zeta
}$ on $\matrixNk$ as follows. Let $A\subset\mathL$ be any
measurable subset, then we define
\[
\mu_{\zeta}\bigl(\{\bM_\lambda\dvtx\lambda\in A\}\bigr)=
\zeta(A).
\]
Let $\zeta$ be a measure on $\mathL$ so that $\mu_\zeta$
satisfies \eqref{eq:finite rates matrix} and let $\mathbf{G}:=\{
(t,G_t)\}\subset\mathbb{R}^+\times\mathL$ be a Poisson point
process with intensity $\mathrm{d}t\otimes\zeta$. Then the construction of the
process $\boldsymbol{\Lambda}$ on $\mathL$ with initial state
$\Lambda_0\in\mathL$ proceeds as in \eqref{eq:PPP formal}
where, for every atom time $t>0$ of $\mathbf{G}$, we define $M_t:=\bM
_{G_t}$. Note that the exchangeability condition from Theorem~\ref
{cor:exchangeable lipschitz} (in coordination with de Finetti) implies
that if the process $\boldsymbol{\Lambda}$ constructed from $\mathbf{G}$
is exchangeable, then the measure $\zeta$ directing $\mathbf{G}$ must
coincide with $\zeta_\nu$ for some measure $\nu$ on $\simplexk$.
%
%co4.1 #&#

\begin{cor}
The process $\boldsymbol{\Lambda}$ based on $\mathbf{G}$ and
$\Lambda
_0$ is exchangeable if and only if $\zeta=\zeta_\nu$, for some
measure $\nu$ on $\simplexk$, and $\Lambda_0$ is an exchangeable
$k$-partition of $\mathbb{N}$.
\end{cor}

By regarding the elements of $\mathbb{N}$ as labeled balls and the
classes of $\lambda\in\mathL$ as labeled boxes, the action $\bM
_\lambda\lambda'$ can be interpreted as the reassignment of each of
the balls to a new class via a cyclic shift by one less than their
class assignment in $\lambda$ (modulo $k$).
The commutative property of this class of maps also implies that the
collection $\{\bM_\lambda\dvtx\lambda\in\mathL\}\subset
\matrixNk
$ is a special subspace of $\matrixNk$. For instance, for every
$\lambda\in\mathL$, $\bM_{\lambda}\bM_{\lambda}^T=I_{\infty
}^k$, where $M^T\in\matrixNk$ denotes the usual matrix transpose of $M$.
\end{example}

%s4.4 #&#
\subsection{Associated \texorpdfstring{$\Delta_k$}{Deltak}-valued Markov process}\label
{section:measure-valued}

We define the \emph{asymptotic frequency} of any $A\subseteq\mathbb
{N}$ by the limit
%
%e4.11 #&#
%
\begin{equation}
\label{eq:asymptotic frequency} |A|:=\lim_{n\rightarrow\infty}\frac{\#(A\cap[n])}{n}, \qquad\mbox
{if it exists}.
\end{equation}
Furthermore, we say $\lambda\in\mathL$ possesses asymptotic
frequency $|\lambda|:=(|\lambda_j|, 1\leq j\leq k)\in\simplexk$,
provided $|\lambda_j|$ exists for every $j=1,\ldots,k$. By de
Finetti's theorem, any exchangeable $k$-partition of $\mathbb{N}$
possesses an asymptotic frequency almost surely. In particular, for
$s\in\simplexk$, $\Lambda\sim\zeta_s$ has $|\Lambda|=s$ with
probability one.

Given a process $\boldsymbol{\Lambda}$ on $\mathL$, its associated
$\simplexk$-valued process is defined by $|\boldsymbol{\Lambda
}|:=(|\Lambda_t|, t\geq0)$, provided $|\Lambda_t|$ exists for all
$t\geq0$ simultaneously. In this section, we show that the associated
$\simplexk$-valued process of any exchangeable Lipschitz partition
process $\boldsymbol{\Lambda}$ exists almost surely and is a Feller process.

Let $\mu$ be the directing measure of an exchangeable Lipschitz
partition process $\boldsymbol{\Lambda}$ and let $\mathbf{M}:=\{
(t,M_t)\}
$ be a Poisson point process with intensity $\mathrm{d}t\otimes\mu$. For $M\in
\matrixNk$, we define the asymptotic frequency of any $M\in\matrixNk
$ as the (column) stochastic matrix $S:=|M|_k$ with $(i,j)$-entry
$S_{ij}:=|M_{ij}|$, provided $|M_{ij}|$ exists for all $i,j=1,\ldots
,k$. We have the following lemmas.
%
%le4.3 #&#

\begin{lemma}\label{prop:atom times}
For every atom time $t>0$ of $\mathbf{M}$, $|M_t|_k$ exists almost surely.
\end{lemma}

\begin{pf}
This is a consequence of Theorem~\ref{cor:exchangeable lipschitz}, by
which, for any atom time $t>0$ of $\mathbf{M}$, each column of $M_t$
is an exchangeable $k$-partition. By de Finetti's theorem, the
asymptotic frequency of each column of $M_t$ exists almost surely.
Since $k<\infty$, $|M_t|_k$ exists a.s.
\end{pf}

For each atom time $t$ of $\mathbf{M}$, we write $S_t:=|M_t|_k$. We
also augment the map $|\cdot|_k$ on $\matrixNk$ by including the
cemetery state $\partial$ in the codomain of $|\cdot|_k$ and defining
$|M|_k=\partial$ if $|M|_k$ does not exist. This makes $|\cdot
|_k\dvtx\matrixNk\rightarrow\stochk\cup\{\partial\}$ a measurable map,
where $\stochk$ is the space of $k\times k$ column stochastic
matrices, that is, $S=(S_{ij}, 1\leq i,j\leq k)\in\stochk$ satisfies
$S_{ij}\geq0$ and $S_{1j}+\cdots+S_{kj}=1$ for all $j=1,\ldots,k$.
%
%le4.4 #&#

\begin{lemma}\label{lemma:cor}
The image $\mathbf{S}:=\{(t,S_t)\}\subseteq\mathbb{R}^+\times
\stochk$ of $\mathbf{M}:=\{(t,M_t)\}\subset\mathbb{R}^+\times
\matrixNk$ by $|\cdot|_k$, that is, $S_t:=|M_t|_k$ for all atom times
$t>0$ of $\mathbf{M}$, is almost surely a Poisson point process with
intensity $\mathrm{d}t\otimes|\mu|_k$, where $|\mu|_k$ denotes the image
measure of $\mu$ by $|\cdot|_k$.
\end{lemma}

\begin{pf}
Let $J\subset[0,\infty)$ denote the subset of atom times of $\mathbf
{M}$. By condition \eqref{eq:finite rates}, $J$ is at most countable
almost surely. By Lemma~\ref{prop:atom times}, $|M_{t}|_k$ exists $\mu
$-almost everywhere for every $t\in J$. Therefore,
\[
\mu \biggl(\bigcup_{t\in J}\bigl \{|M_t|_k=
\partial\bigr \} \biggr)\leq\sum_{t\in J}\mu \bigl(
\bigl \{|M_t|_k=\partial\bigr \} \bigr)=0,
\]
and $\mathbf{S}$ is almost surely a subset of $\mathbb{R}^+\times
\stochk$. That $\mathbf{S}$ is a Poisson point process with the
appropriate intensity is clear as it is the image of the Poisson point
process $\mathbf{M}$ by the measurable function $|\cdot|_k$.
\end{pf}

%th4.2 #&#
%
\begin{them}\label{thm:measure-valued}
The associated $\simplexk$-valued process $|\boldsymbol{\Lambda
}|:=(|\Lambda_t|, t\geq0)$ of an exchangeable Lipschitz partition
process $\boldsymbol{\Lambda}$ exists almost surely and is a Feller
process on $\simplexk$.
\end{them}
\begin{pf}
By exchangeability of $\boldsymbol{\Lambda}$, the asymptotic frequency
$|\Lambda_t|$ exists for all fixed times $t>0$, with probability one.
In order for $|\boldsymbol{\Lambda}|$ to exist on $\simplexk$, we must
show that, with probability one, $|\Lambda_t|$ exists for all $t>0$
simultaneously.
Let $\mathbf{M}$ be a Poisson point process that determines the jumps
of $\boldsymbol{\Lambda}$ by \eqref{eq:PPP formal}. For each $n\in
\mathbb{N}$, let $\mathcal{D}_n$ be the dyadic rationals in $[0,1]$.
Then $|\Lambda_t|$ exists almost surely on $\bigcup_{n\in\mathbb
{N}}\mathcal{D}_n$, which is dense in $[0,1]$. Existence of
$|\boldsymbol
{\Lambda}|$ now follows by density, c\`adl\`ag paths of $\boldsymbol
{\Lambda}$, the Poisson point process construction of $\boldsymbol
{\Lambda}$ via $\mathbf{M}$, and Lemmas~\ref{prop:atom times} and~\ref{lemma:cor}.

The Feller property follows from Corollary~\ref{cor:lip-ppp} and Lemma~\ref{lemma:cor}. As in the general case, in which $\phi^{\infty}$ is
a Feller process on $\Lip(\mathL)$, we can construct a Feller
process $\mathbf{Q}:=(Q_t, t\geq0)$ on $\matrixNk$ such that
$\mathbf{Q}^{\sigma}:=(Q^{\sigma}_t, t\geq0)\equalinlaw\mathbf
{Q}$, for all $\sigma\dvtx\mathbb{N}\rightarrow\mathbb{N}$ fixing all
but finitely many $n\in\mathbb{N}$. By Corollary~\ref{cor:lip-ppp},
the semigroup of $\boldsymbol{\Lambda}$ satisfies $\mathbf
{P}_tg(\lambda
):=\mathbb{E}_{\lambda}g(Q_t\lambda)$. Furthermore, by Lemma~\ref
{lemma:cor} and the argument to show that $|\boldsymbol{\Lambda}|$
exists, the projection $|\mathbf{Q}|:=(|Q_t|_k, t\geq0)$ into
$\stochk$ exists almost surely and $|\boldsymbol{\Lambda}|$ satisfies
$\Lambda_t:=|Q_t|_k|\Lambda_0|$ for all $t>0$. The Feller property is
a consequence of Lipschitz continuity of the linear map $S\dvtx
\simplexk
\rightarrow\simplexk$ determined by any $S\in\stochk$.\end{pf}

%re4.3 #&#

\begin{rmk}
A detailed proof of Theorem~\ref{thm:measure-valued} is technical and
provides no new insights.
Essentially, existence of $|\boldsymbol{\Lambda}|$ is a consequence of
regularity of the paths of $\boldsymbol{\Lambda}$ and density of the
countable set of dyadic rationals. The Feller property follows by
Lipschitz continuity of maps determined by stochastic matrices. For a
blueprint of the proof, we point the reader to \cite{CraneLalley2012b}.
\end{rmk}

%s4.5 #&#
\subsection{A non-exchangeable Lipschitz process}\label
{section:non-exchangeable process}

The processes in the above examples are exchangeable Lipschitz
partition processes.
We now show an example of a Lipschitz partition process that is not
exchangeable, and whose directing measure is not confined to the
subspace $\matrixNk$.

Let $\mathbf{A}:=(M^i_j, i\in[k], j\geq0)$ be an array of elements
in $\mathL$. Given $\mathbf{A}$, we define $F_{\mathbf{A}}=F\in
\Phi$ by $F(\lambda)=\mathbf{A}_\lambda\lambda$, where $\mathbf
{A}_\lambda\in\matrixNk$ is defined as follows. For every $i\in
[k]$, we put
%
%e4.12 #&#
%
\begin{equation}
\label{eq:minima}m_i:=\cases{ %
%\begin{array}{cc}
\min\{n\in
\mathbb{N}\dvtx n\in\lambda_i\},&\quad$\lambda_i\neq
\emptyset$,
\cr
0,&\quad$\mbox{otherwise}$. %\end{array}
%%
%\right.
}
\end{equation}
For each $i=1,\ldots,k$, we put $\mathbf{A}^i_\lambda=M^i_{m_i}$ and
let $\mathbf{A}_\lambda:=(\mathbf{A}^1_\lambda,\ldots,\mathbf
{A}^k_\lambda)\in\matrixNk$. It should be clear that, as specified,
$F$ need not be \emph{strongly} Lipschitz.
%
%pr4.3 #&#

\begin{prop}\label{prop:Lipschitz map}
The map $F_{\mathbf{A}}$ defined above is Lipschitz continuous.
\end{prop}

\begin{pf}
Take any $\lambda,\lambda'\in\mathL$ and let $r=-\log_2
d(\lambda,\lambda')\in\mathbb{N}\cup\{\infty\}$. Let
$0<m_{(1)}<m_{(2)}<\cdots<m_{(k')}\leq r$ and
$0<m_{(1)}'<m_{(2)}'<\cdots<m_{(k'')}'\leq r$ be the minima \eqref
{eq:minima} of $\lambda$ and $\lambda'$ (respectively) that are
greater than zero but not greater than $r$. Since $I_r^k\lambda
=I_r^k\lambda'$ by definition \eqref{eq:ultrametric}, we must have
$k''=k'$ and $m_{(i)}=m_{(i)}'$ for all $1\leq i\leq k'$. It follows
that $\mathbf{A}_\lambda I_r^k=\mathbf{A}_{\lambda'}I_r^k$ and
\begin{eqnarray*}
F_{\mathbf{A}}(\lambda)_{|[r]}&=&\bigl(I_r^k
\mathbf{A}_\lambda\bigr)\lambda =\bigl(\mathbf{A}_\lambda
I_r^k\bigr)\lambda=\mathbf{A}_\lambda
\bigl(I_r^k\lambda \bigr)=\mathbf{A}_\lambda
\bigl(I_r^k\lambda'\bigr)=\bigl(
\mathbf{A}_\lambda I_r^k\bigr)
\lambda'=\bigl(\mathbf{A}_{\lambda'} I_r^k
\bigr)\lambda'\\
&=&F_{\mathbf
{A}}\bigl(\lambda'
\bigr)_{|[r]}.
\end{eqnarray*}
As this must hold for all $\lambda,\lambda'\in\mathL$, it
follows that $F_{\mathbf{A}}$ is Lipschitz continuous.
\end{pf}

Now, we construct a measure on $\Lip(\mathL)$ using the above
observation. In particular, for every $j\geq0$, let $\nu_{j}$ be a
measure on $\simplexk$ such that
%
%e4.13 #&#
%
\begin{equation}
\label{eq:finite lambda rates}\zeta^{(n)}_{\nu
_{j}}\bigl({\mathcal{L}_{[n]:k}}
\setminus\bigl\{E_i^{(n)}\bigr\} \bigr)<\infty\qquad
\mbox{for all }i=1,\ldots,k, \mbox{ for all }n\in \mathbb{N},
\end{equation}
where $E_i^{(n)}\in{\mathcal{L}_{[n]:k}} $ is the
$k$-partition of $[n]$ with all elements labeled $i$, as in the proof
of Proposition~\ref{cor:isomorphism}. We define the measure $\mu$ on
$k\times\infty$ arrays of independent $k$-partitions (as $\mathbf
{A}$ above) for which the partition in the $i$th row and $j$th column
has distribution $\zeta_{\nu_{j}}$. We then define $\varphi_\mu$ as
the measure on $\Lip(\mathL)$ induced by the random array
$\mathbf{A}$ with distribution $\mu$ and the map $F\in\Lip
(\mathL)$ associated to $\mathbf{A}$ by the above discussion. We
let $\mathbf{F}$ be a Poisson point process with intensity $\mathrm{d}t\otimes
\varphi_\mu$ and construct $\boldsymbol{\Lambda}$ on $\mathL$ as
in \eqref{eq:PPP formal}.

In the following proposition, let $\mathbf{F}:=\{(t,F^t)\}$ be a
Poisson point process with intensity $\mathrm{d}t\otimes\varphi_\mu$, which
is determined by a Poisson point process $\mathbf{A}:=\{(t,A_t)\}$
with intensity $\mathrm{d}t\otimes\mu$, where each $A_t$ is a random $k\times
\infty$ array. In particular, for each atom time $t>0$ of $\mathbf
{F}$, we put $F^t:=F_{A_t}$, as defined above.
%
%pr4.4 #&#

\begin{prop}\label{prop:non-exchangeable}
$\boldsymbol{\Lambda}$ constructed from $\mathbf{F}$ is a Feller process
on $\mathL$. If, in addition,
$\nu_i\neq\nu_j$ for some $1\leq i<j<\infty$, then $\boldsymbol
{\Lambda}$ is not exchangeable.
\end{prop}

\begin{pf}
For every $n\in\mathbb{N}$ and atom time $t>0$ of $\mathbf{A}$, the
restriction $\Lambda_{t|[n]}$ depends only on the first $n+1$ columns
of any $A_t$. By assumption \eqref{eq:finite lambda rates} on the
underlying directing measures $\zeta_{\nu_{j}}$, $\varphi_\mu$
satisfies \eqref{eq:finite rates}. Theorem~\ref{thm:feller} and
Proposition~\ref{prop:Lipschitz map} now imply that $\boldsymbol
{\Lambda
}$ is a Feller process.

Non-exchangeability of $\boldsymbol{\Lambda}$ under the stated condition
is clear: since $\nu_i\neq\nu_j$ implies $\zeta_{\nu_i}^{(n)}\neq
\zeta_{\nu_j}^{(n)}$ for all $n\in\mathbb{N}$, then the jump rates
from a state with $m_{i'}=i$ and $m_{j'}=j$ differ from the jump rates
from a state with $m_{i'}=j$ and $m_{j'}=i$, for any $1\leq i'\neq
j'\leq k$.
\end{pf}

%s5 #&#
\section{Discrete-time processes}\label{section:discrete-time processes}

From our previous discussion of continuous-time processes, we need not
prove anything further for discrete-time chains; but we make some
observations specific to the discrete-time case. Throughout this
section, all measures on $\Lip(\mathL)$ and/or $\matrixNk$ are
\emph{probability measures}.

First, given a probability measure $\varphi$ on $\Lip(\mathL)$,
we construct a Markov chain $\boldsymbol{\Lambda}:=(\Lambda_m, m\geq0)$
with initial state $\Lambda_0\in\mathL$ by taking
$F_1,F_2,\ldots$ i.i.d. with law $\varphi$ and defining
%
%e5.1 #&#
%
\begin{equation}
\label{eq:discrete-time formal}\Lambda_m=F_m(\Lambda
_{m-1})=(F_m\circ F_{m-1}\circ\cdots\circ
F_1) (\Lambda_0),\qquad\mbox{for each }m\geq1.
\end{equation}
Constructed this way, $\boldsymbol{\Lambda}$ is a Markov chain on
$\mathL$. Furthermore, by Lipschitz continuity of the maps
$F_1,F_2,\ldots$\,, the finite restrictions $(\boldsymbol{\Lambda}_{|[n]},
n\in\mathbb{N})$ are finite state space Markov chains. The following
corollary follows from arguments in the continuous-time case.

%co5.1 #&#

\begin{cor}\label{cor:discrete-time Markov}
Let $\boldsymbol{\Lambda}$ constructed in \eqref{eq:discrete-time
formal} be exchangeable. Then we have the following.
\begin{itemize}
\item$\varphi$ is supported on $\mathcal{F}\cap\Sigma$ and we can
assume, without loss of generality, that $\varphi$ is a probability
measure on $\matrixNk$.
\item The $\simplexk$-valued Markov chain $|\boldsymbol{\Lambda
}|=(|\Lambda_m|, m\geq0)$ exists almost surely and can be constructed
as in \eqref{eq:discrete-time formal} from $S_1,S_2,\ldots$ i.i.d.
$|\varphi|_k$, the measure induced by $\varphi$ on $\stochk$ through
the map $|\cdot|_k$. In particular, $|\boldsymbol{\Lambda
}|\equalinlaw
\mathbf{D}:=(D_m, m\geq0)$, where
\[
D_m:=S_m\cdots S_1 D_0,\qquad
m\geq1,
\]
for $D_0:=|\Lambda_0|$ and $S_1,S_2,\ldots$ i.i.d. $|\varphi|_k$.
\end{itemize}
\end{cor}

%s6 #&#
\section{Concluding remarks}\label{section:concluding remarks}

To conclude, we remark about the projection of Lipschitz partition
processes into $\partitionsNk$ and discuss more general aspects of
partition-valued processes.

%s6.1 #&#
\subsection{Associated Lipschitz partition processes on $\partitionsNk
$}\label{section:partition-processes}

Let $\varphi$ be the directing measure of a Lipschitz partition
process $\boldsymbol{\Lambda}$ on $\mathL$. Intuitively, the
projection $\mathcal{B}_{\infty}(\boldsymbol{\Lambda}):=(\mathcal
{B}_{\infty}(\Lambda_t), t\geq0)$ into $\partitionsNk$ is, itself,
a Markov process as long as $\varphi$ treats the classes of every
$\lambda\in\mathL$ ``symmetrically.'' In particular, for any
permutation $\gamma\dvtx[k]\rightarrow[k]$, let us define $\Gamma
\in
\matrixNk$ as the $k\times k$ partition operator with entries
\[
\Gamma_{ij}=\cases{ %
%\begin{array}{cc}
\mathbb{N}, & \quad$
\gamma(i)=j$,
\cr
\emptyset, & \quad$\mbox{otherwise}$. %\end{array}
%%
%\right.
}
\]
The matrix $\Gamma$ acts on $\mathL$ by relabeling classes; that
is, for any $\lambda:=(\lambda_i, 1\leq i\leq k)\in\mathL$,
$\Gamma\lambda:=(\lambda_{\gamma(i)}, 1\leq i\leq k)$. Therefore,
the projection $\boldsymbol{\Pi}:=\mathcal{B}_{\infty}(\boldsymbol
{\Lambda
})$ into $\partitionsNk$ is a Markov process if and only if, for every
$\lambda\in\mathL$ and every measurable subset $C\subseteq
\partitionsNk$, $\varphi$ assigns equal measure to the events $\{F\in
\Phi\dvtx F(\lambda)\in\mathcal{B}_{\infty}^{-1}(C)\}$ and $\{
F\in
\Phi\dvtx F(\Gamma\lambda)\in\mathcal{B}_{\infty}^{-1}(C)\}$,
for all
$\gamma\in\symmetrick$. Moreover, if $\boldsymbol{\Pi}$ is a Markov
process, then it fulfills the Feller property.

By the preceding discussion, we can generate a Lipschitz partition
process on $\partitionsNk$ by projecting a process $\boldsymbol
{\Lambda
}$ that treats labels symmetrically. The projection $\mathcal
{B}_{\infty}(\boldsymbol{\Lambda})$ is a Feller process; and, if
$\boldsymbol{\Lambda}$ is exchangeable, then so is $\mathcal
{B}_{\infty
}(\boldsymbol{\Lambda})$.

%s6.2 #&#
\subsection{Existence and related notions}

Sections~\ref{section:examples} and~\ref{section:non-exchangeable
process} contain explicit examples of exchangeable and non-exchangeable
Lipschitz partition processes.
These examples confirm that Lipschitz partition processes exist, and
their Poisson point process construction lends insight into their
behavior. The Poisson point process construction is also useful in
simulation and Markov chain Monte Carlo sampling.

There remain broader questions surrounding existence of measures
satisfying \eqref{eq:finite rates matrix}, as well as more general
partition-valued Markov processes. We undertake some of these questions
elsewhere:
we characterize exchangeable Feller processes on $\partitionsNk$ in
\cite{Crane2013CP};
we show the cutoff phenomenon for a class of these chains in \cite
{CraneLalley2012a}; and we study exchangeable processes without the
Feller property in
\cite{CraneLalley2012b}.

%\begin{appendix}
%\section{}
%\end{appendix}

\section*{Acknowledgements}
The author is partially supported by NSF grant DMS-1308899 and NSA
grant H98230-13-1-0299.

%\begin{supplement}%[id=suppA]
%\slink[doi]{10.3150/00-BEJXXXXSUPP} %[doi,text={...}] - jei reikia
%suskaldyti doi
%\sfilename{BEJ000\_supp.pdf}
%\end{supplement}

%\bibitem[\protect\citeauthoryear{}{()}]{r1}
%\bibitem{r1}
% imsref loaded by audrone.aklyte, 2014-04-07 09:15:10
% imsref loaded by audrone.aklyte, 2014-04-07 09:18:33
%

\printhistory

\end{document}